\documentclass{article}
\usepackage {amsmath, amsfonts, amssymb, mathrsfs, amsthm,stmaryrd, sectsty,hyphenat,etoolbox,tabu,calc,graphicx,url,enumitem}
\usepackage{palatino}
\usepackage[all]{xy}
\usepackage[margin=1.5in]{geometry}

\def\mf#1{\mathfrak{#1}}
\def\mc#1{\mathcal{#1}}
\def\mb#1{\mathbb{#1}}
\def\tx#1{\textrm{#1}}

\def\ms#1{\mathsf{#1}}

\def\R{\mathbb{R}}

\def\C{\mathbb{C}}
\def\Q{\mathbb{Q}}

\def\Z{\mathbb{Z}}

\def\F{\mathbb{F}}
\def\lmod{\setminus}

\def\ol#1{\overline{#1}}

\def\hat{\widehat}

\def\lrw{\longrightarrow}
\def\llw{\longleftarrow}

\def\<{\langle}
\def\>{\rangle}

\newenvironment{mytitle}
{\begin{center}\large\sc}
{\end{center}}

\newtheorem{thm}{Theorem}[section]

\newtheorem{cor}[thm]{Corollary}

\newtheorem{cnj}[thm]{Conjecture}

\theoremstyle{definition}
\newtheorem{rem}[thm]{Remark}

\newtheorem{qst}[thm]{Question}

\numberwithin{equation}{section}
\sectionfont{\center\sc\normalsize}
\subsectionfont{\bf\normalsize}

\AtBeginDocument{%
   \def\MR#1{}
}

\begin{document}

\begin{mytitle}
Representations of reductive groups over local fields\\
\end{mytitle}

\begin{center}
Tasho Kaletha
\end{center}

\begin{abstract}
We discuss progress towards the classification of irreducible admissible representations of reductive groups over non-archimedean local fields and the local Langlands correspondence. We also state some (partly conjectural) compatibility properties of the refined local Langlands correspondence.
\end{abstract}

\let\thefootnote\relax
\footnotetext{This work was partially supported by NSF Grant DMS-1801687 and a Simons Fellowship.}

\setcounter{secnumdepth}{3}
\setcounter{tocdepth}{2}
\tableofcontents

\section*{Introduction}

Let $F$ be a local field, i.e. a finite extension of the field $\R$ of real numbers, or the field $\Q_p$ of $p$-adic numbers, or the field $\F_p((t))$ of Laurent series over a finite field. Let $G$ be a connected reductive $F$-group. Motivated by the theory of automorphic forms, the study of irreducible admissible representations of the topological group $G(F)$ with complex coefficients has been an active area of research since the pioneering work of Bargmann on $\tx{SL}_2(\R)$. Some of the main problems in this area are
\begin{enumerate}[itemsep=0.5ex,topsep=0.6ex]
	\item The classification of irreducible admissible representations.
	\item The determination of their character functions.
\end{enumerate}
In addition, motivated by Langlands' conjectures, one can add
\begin{enumerate}[itemsep=0.5ex,topsep=0.6ex,resume]
	\item The relation with representations of the Galois/Weil group of $F$.
	\item The proof of character identities stemming from endoscopy and more general functoriality.
	\item The appropriate normalization of intertwining operators for parabolic induction.
\end{enumerate}
In this note I would like to discuss progress towards some of these questions. There are of course many other interesting questions, such as for example the classification of unitary representations, about which I will not say anything.\\[-6pt]

\textbf{Acknowledgements:} We thank Michael Harris, Robert Kottwitz, and Gopal Prasad, for helpful remarks during the preparation of this article.

\section{Classification of irreducible representations and characters}

\subsection{The archimedean case}

The archimedean case, where $F$ is a finite extension of $\R$, thus equal to $\R$ or $\C$, has been largely resolved by the work of Harish-Chandra, Langlands, Shelstad, and others. I will discuss it briefly, because it will serve as a useful guide to the non-archimedean case. 

One of Harish-Chandra's fundamental contributions was the introduction of the notion of a discrete series representation, i.e. those unitary representations whose matrix coefficients are square\-integrable modulo center, and the classification of such representations. Slightly more generally, one considers essentially discrete series representations -- those which become discrete series after tensor product with a character of $G(F)$. His theorem can roughly be stated as follows (this is a reformulation due to Langlands \cite{Lan89} and is easily seen to be equivalent to the original formulation).

\begin{thm} \label{thm:hc1}
The set of isomorphism classes of essentially discrete series representations of $G(F)$ is in a natural bijection with the set of $G(F)$-conjugacy classes of triples $(S,B,\theta)$, where $S \subset G$ is an elliptic (i.e. anisotropic modulo center) maximal torus, $B$ is a Borel subgroup of $G_{\bar F}$ containing $S$, and $\theta$ is a character of $S(F)$ whose differential is $B$-dominant.
\end{thm}

The group $G$ may fail to have an elliptic maximal torus. For example, when $F=\C$, such a torus never exists, unless $G$ itself is a torus. Even when $F=\R$ an elliptic maximal torus may not exist, as for example in the case of $\tx{SL}_n$ for $n>2$. When an elliptic maximal torus $S \subset G$ exists, it is unique up to $G(\R)$-conjugacy. The corresponding maximal torus $S_\tx{sc}$ of the simply connected cover of the derived subgroup of $G$ is anisotropic, and the restriction of $\theta$ to it is an algebraic character, i.e. an element of the coweight lattice $X^*(S_\tx{sc})$ of the absolute root system of $G$ relative to $S$, so it makes sense to ask that $d\theta$ be $B$-dominant. Moreover, when this differential is a regular element of the weight lattice, it uniquely determines $B$, so the classifying datum is just a $G(\R)$-conjugacy class of pairs $(S,\theta)$. We might be tempted to call the corresponding essentially discrete series representation ``regular'', a notion that will find its analog in the non\-archimedean case.

The essentially discrete series representation $\pi_{(S,B,\theta)}$ associated to the triple $(S,B,\theta)$ by the above theorem can be specified by its character function. This uses another fundamental result of Harish-Chandra (valid for an arbitrary local field $F$ of characteristic zero, and extended to local fields of positive characteristic by \cite{CGH14} under some assumptions): the fact that the character distribution
\[ f \mapsto \tx{tr}\,\pi(f),\qquad \pi(f)v = \int_{G(F)} f(g)\pi(g)vdg \]
of an admissible representation $\pi$ of $G(F)$ is representable by a locally integrable function $\Theta_\pi : G(F) \to \C$. Just like in the case of finite groups, this function determines $\pi$ up to equivalence. In fact, when $F=\R$ and $\pi$ is essentially discrete series, already the restriction of $\Theta_\pi$ to $S(\R)$ determines $\pi$. More precisely:

\begin{thm} \label{thm:hc2}
$\pi_{(S,B,\theta)}$ is the unique essentially discrete series representations such that for all $s \in S(F) \cap G(F)_\tx{reg}$
\[ \Theta_{\pi_{(S,B,\theta)}}(s) = (-1)^{q(G)}\sum_{w \in N(S,G)(\R)/S(\R)} \frac{\theta(s^w)}{\prod\limits_{\alpha >0} (1-\alpha(s^w)^{-1})}, \]
where $q(G)$ is half of the dimension of the symmetric space of $G(\R)$ and $\alpha>0$ indicates the product over all those absolute roots with respect to $S$ that are positive with respect to $B$.	
\end{thm}
 
The classification of (essentially) discrete series representations of $G(F)$ is a key step in the classification of all irreducible admissible representations of $G(F)$. The next step is the classification of the irreducible tempered representations, i.e. those unitary representations whose matrix coefficients are almost square\-integrable modulo center. Harish-Chandra has shown that these are precisely the irreducible constituents of parabolically induced discrete series representations of Levi subgroups of $G$. Moreover, the theory of the $R$-group due again to Harish-Chandra provides a description of the various irreducible constituents of such a parabolic induction, and hence a full classification of the tempered representations, provided one has suitably normalized the intertwining operators. What the right normalization is has been conjectured by Langlands for any local field. For archimedean local fields Arthur \cite[\S3]{ArtIOR1} proved Langlands' conjecture, while for non\-archimedean local fields of characteristic zero Arthur \cite[\S4]{ArtIOR1} proved abstractly that a normalization exists, without being able to prove that it is provided by Langlands' formula.

The final step is the Langlands classification theorem, which states that every irreducible admissible representation is equivalent to one of the form $j_P^G(\sigma\otimes\nu)$. Here $P$ is a parabolic subgroup of $G$; we denote by $M$ the Levi quotient of $P$ and by $A_M$ the maximal split torus in the center of $M$; $\sigma$ is a tempered representation of $M(F)$ and $\nu$ is an element of $X^*(A_M) \otimes \R$ that lies in the acute open cone associated to $P$, and which is identified with a character of $M(F)$ using the exponential map; and $j_P^G$ is the unique irreducible quotient of the parabolic induction $i_P^G(\sigma\otimes\nu)$. It is known that two such representations are equivalent if and only if their Langlands data $(P,\sigma,\nu)$ are $G(F)$-conjugate.

\subsection{The non-archimedean case}

Consider now a finite extension $F$ of $\Q_p$ or $\F_p((t))$. Harish-Chandra's classification of tempered representations in terms of discrete series representations of Levi subgroups by means of the theory of the $R$-group, and Langlands' classification theorem, continue to hold, with minor modifications to their statements and proofs, cf. \cite[CH. VII]{Ren10}. This is an instance of Harish-Chandra's ``Lefschetz principle'', which is the philosophy that the representation theory of real and $p$-adic groups (and even the automorphic representations of adele groups) exhibit parallel behavior, despite the stark differences in the fine structure of these groups. But, unlike in the archimedean case, our understanding of the discrete series representations in the non\-archimedean case is less developed, and the classification of these representations is at this moment incomplete.

A special subclass of the discrete series is made out of the supercuspidal representations, which are those whose matrix coefficients have compact support modulo the center. Real reductive groups do not have such representations, except for the trivial case of tori. The last 30 years have seen a significant improvement of our understanding of supercuspidal representations, beginning with the work of Moy--Prasad \cite{MP94,MP96} and Morris \cite{Mor88,Mor89} in the case of depth zero, the constructions of general depth supercuspidal representations due to Adler \cite{Ad98} and Yu \cite{Yu01}, and the exhaustion results of Kim \cite{Kim07} and Fintzen \cite{Fin21}. These works rely crucially on the filtrations of the topological group $G(F)$ coming from Bruhat--Tits theory and its extensions by Moy--Prasad (an example is the filtration of the compact group $\tx{SL}_2(\Z_p)$ by congruence subgroups), and as such are very much a $p$-adic phenomenon with no clear analog in the archimedean case. 

The most comprehensive construction, due to Yu, produces supercuspidal representations out of what is nowadays customarily called ``Yu-data'', rather complicated structures consisting among other things of a tower of twisted Levi subgroups, a depth\-zero supercuspidal representation of the smallest subgroup, and a sequence of characters of each subgroup subject to a genericity condition. Fintzen's result shows that all supercuspidal representations arise from this construction when $G$ is tamely ramified and $p$ does not divide the order of the absolute Weyl group of $G$. Work of Hakim--Murnaghan \cite{HM08} defines an explicit equivalence relation on the set of Yu-data which describes when two data produce the same representation. These results amount to a classification of all supercuspidal representations of $G(F)$ under the given conditions on $G$, in terms of equivalence classes of Yu-data. However, a simpler classification may be desirable. In fact, with Harish-Chandra's Lefschetz principle and his work for real groups in mind, we would ideally like a classification in terms of objects close to $G(F)$-conjugacy classes of pairs $(S,\theta)$ consisting of an elliptic maximal torus $S \subset G$ and a character $\theta$ of it.

Moy--Prasad introduce the notion of \emph{depth} of a representation and show that an irreducible depth\-zero supercuspidal representation always arises via compact induction from a maximal open and compact\-mod\-center subgroup $G(F)_x$ -- the stabilizer of a vertex in the Bruhat--Tits building -- of an irreducible representation $\sigma$ of $G(F)_x$ with the following very special property: $G(F)_x$ has a natural quotient that is the group $\ms{G}_x(k_F)$ of $k_F$-points of a usually disconnected reductive $k_F$-group $\ms{G}_x$, and $\sigma$ is required to factor through this quotient and moreover the restriction to the identity component $\ms{G}_x^\circ(k_F)$ must contain a cuspidal representation of this finite group of Lie type; here $k_F$ is the residue field of $F$. In this way, the representation theory of finite groups of Lie type (including disconnected ones) is reflected in the representation theory of reductive $p$-adic groups. Given a connected reductive $k_F$-group $\ms{G}$ and a pair $(\ms{S},\theta)$ of a maximal torus $\ms{S} \subset \ms{G}$ and a character $\theta$ of $\ms{S}(k_F)$ (customarily taking $\ell$-adic values), the construction of Deligne--Lusztig \cite{DL76} assigns a virtual representation $R_{\ms{S},\theta}$ of $\ms{G}(k_F)$. In general this virtual representation is not an actual representation (even up to sign), but quite often it is. More precisely, Deligne--Lusztig define the notions of a character $\theta$ to be ``non-singular'' and in ``general position'', which are dual to the notions of a semi\-simple element in a connected reductive group to be ``regular'' and ``strongly regular''. They show \cite[Remark 9.15.1]{DL76} that $R_{\ms{S},\theta}$ is an actual representation (up to a well\-understood sign) whenever $\theta$ is non\-singular (originally under a certain affineness assumption, which was later shown to always hold by He \cite{He08}). Moreover, Deligne--Lusztig show that the representation $\pm R_{\ms{S},\theta}$ is irreducible when $\theta$ is in general position, and cuspidal if $\ms{S}$ is elliptic. 

These results are both encouraging for our quest to parameterize supercuspidal representations in terms of pairs $(S,\theta)$, but also cautioning us that there will be supercuspidal representations that do not obey such a parameterization. More precisely, in \cite[\S3]{KalRSP} I define the notion of a ``regular'' supercuspidal representation, which is one that arises from Yu's construction and for which the depth\-zero part of the Yu-datum comes from, via the results of Moy--Prasad and Deligne--Lusztig, a character in general position, and then prove the following classification.

\begin{thm}[{\cite[Cor. 3.7.10]{KalRSP}}] \label{thm:k1}
Assume that $G$ splits over a tame extension of $F$ and $p$ does not divide the order of the Weyl group of $G$. The set of isomorphism classes of regular supercuspidal representations of $G(F)$ is in a natural bijection with the set of $G(F)$-conjugacy classes of pairs $(S,\theta)$, where $S$ is an elliptic maximal torus that splits over a tame extension, and $\theta$ is a regular character of $S(F)$.
\end{thm}

\begin{rem} \label{rem:r1}
\ \\[-15pt]
\begin{enumerate}%
	\item The condition on $p$ can be weakened; I have opted here for the one which is easiest to state.
	\item I have not explicated here the definition of a ``regular'' character $\theta$, but the main point is that it is an explicit Lie\-theoretic condition, essentially amounting to the stabilizer in $N_G(S)(F)/S(F)$ of the restriction $\theta^0$ of $\theta$ to the Iwahori subgroup of $S(F)$ being trivial. For details, cf. \cite[Def. 3.7.5]{KalRSP}.
	\item In the definition of regular supercuspidal representation, ``general position'' should be taken with respect to the $p$-adic group $G$, which is slightly stronger than taking it with respect to the finite group of Lie type $\ms{G}_x^\circ$.
	\item One of the useful properties of this theorem is that it does not reference the fine structure of the topological group $G(F)$ coming from the $p$-adic field $F$, such as the various filtrations coming from Bruhat--Tits and Moy--Prasad theory.
	\item Continuing the previous point, this theorem is in fact rather analogous to Theorem \ref{thm:hc1} restricted to regular discrete series representations (in the sense of the previous subsection). In this way, it establishes the Harish-Chandra Lefschetz principle among a wide class of discrete series representations, setting up a parallel between the regular discrete series representations of real reductive groups and the regular supercuspidal representations of $p$-adic reductive groups.
	\item One important difference between the real and $p$-adic cases is that while in the real case $S$ is unique up to $G(F)$-conjugacy, in the non\-archimedean case there usually are (finitely) many different $G(F)$-conjugacy classes (in fact even isomorphism classes) of elliptic tamely ramified maximal tori of $G$. Moreover, in the non\-archimedean case, elliptic maximal tori always exist, and supercuspidal representations always exist, cf. \cite{Kret12}.
\end{enumerate}
\end{rem}

There is also an analog of Theorem \ref{thm:hc2} which we will discuss in a moment, but before doing so we briefly consider going beyond the ``regular'' case. We can impose on $\theta$ the $p$-adic analog of Deligne--Lusztig's ``non\-singular'' condition, which is weaker than the condition of being regular. The arguments involved in the proof of Theorem \ref{thm:k1} still apply and produce a supercuspidal representation $\pi_{(S,\theta)}$, which may however be reducible, in fact a direct sum of finitely many irreducible supercuspidal representations. 

The irreducible representations obtained this way, i.e. the irreducible constituents of $\pi_{(S,\theta)}$ for all possible pairs $(S,\theta)$, can be characterized in the same way as the regular ones, but where we replace the ``general position'' requirement with a ``non\-singular'' requirement \cite[Def. 3.1.1]{KalSLP}. We can thus call these supercuspidal representations ``non\-singular'' (although a better term might be ``semi\-simple'', as a contrast to the concept of a unipotent supercuspidal representation). It may not be clear at first sight why this class of representations is interesting, beyond it being a generalization of the class of regular supercuspidal representations. The main interest in them comes from the fact that these are precisely those supercuspidal representations whose Langlands parameters are ``supercuspidal'', i.e. discrete with trivial monodromy, at least when $p$ does not divide the order of the Weyl group and according to the construction of \cite{KalSLP}; we will discuss this point in the next section.

The classification of all irreducible non-singular supercuspidal representations reduces, via the analog of Theorem \ref{thm:k1}, to the study of the internal structure of the representations $\pi_{(S,\theta)}$, i.e. its decomposition into irreducible factors and their multiplicities. The situation is made subtle by the fact that an irreducible constituent of $\pi_{(S,\theta)}$ may occur with multiplicity greater than $1$, which is a phenomenon that does not occur for connected reductive groups over finite fields. The study of the internal structure of $\pi_{(S,\theta)}$ reduces to the case of depth\-zero, where it relies on  geometric intertwining operators acting on Deligne--Lusztig induction in the setting of disconnected groups, building on the work of Bonnaf\'e--Dat--Rouquier \cite{BDR17}. These operators must be suitably normalized. Using square brackets to denote sets of irreducible constituents, the following classification result is proved in \cite[\S3]{KalSLP}.

\begin{thm}[{\cite[Fact 2.4.11, Proposition 3.2.4, Proposition 3.4.6, Corollary 3.4.7]{KalSLP}}] \label{thm:k2}
Assume that $G$ splits over a tame extension of $F$ and $p$ does not divide the order of the Weyl group of $G$. 
\begin{enumerate}%
	\item The set of normalizations of the geometric intertwining operators is a non\-empty torsor under the Pontryagin dual of the finite abelian group $N_G(S)(F)_{\theta^0}/S(F)$ (cf. Remark \ref{rem:r1}(2))
	\item Any such normalization provides a multiplicity\-preserving bijection 
	\[ [\pi_{(S,\theta)}] \leftrightarrow [\tx{Irr}_\theta(N_G(S)(F)_\theta], \]
	where on the right we have those irreducible representations of $N_G(S)(F)_\theta$ whose restriction to $S(F)$ is $\theta$-isotypic.
\end{enumerate}
\end{thm}
The existence of normalized intertwining operators, which is part of the first point, is formally analogous to Arthur's result \cite[\S4]{ArtIOR1} on the existence of suitable normalization of standard intertwining operators between parabolically induced representations. As in Arthur's situation, I have not been able to provide a specific normalization. In fact, at the moment there is not even a conjectural expectation of what a good normalization might look like. The fact that the decomposition of the supercuspidal representation $\pi_{(S,\theta)}$ is formally analogous to the decomposition of a parabolic induction (via standard intertwining operators and the $R$-group), is quite intriguing.

The above theorem implies that the set of isomorphism classes of non\-singular supercuspidal representations of $G(F)$ is in bijection with the set of $G(F)$-conjugacy classes of triples $(S,\theta,\rho)$, where $S$ is an elliptic maximal torus that splits over a tame extension, $\theta$ is a non\-singular character of $S(F)$, and $\rho$ is an irreducible representation of $N_G(S)(F)_\theta$ whose restriction to $S(F)$ is $\theta$-isotypic. The bijection is at the moment not completely natural, due the lack of natural normalization of the intertwining operators.

We now come to the non\-archimedean analog of Theorem \ref{thm:hc2}. It is based on work of Adler, DeBacker, Reeder, and Spice, \cite{AS09}, \cite{DR09}, \cite{DS18}, \cite{Spice18}, \cite{Spice21}, and is ultimately formulated in \cite{FKS}. First, we state a simpler version.

\begin{thm}[{\cite[Proposition 4.3.2]{FKS}}] \label{thm:k3}
Let $\pi_{(S,\theta)}$ be the (possibly reducible) supercuspidal representation associated to a pair $(S,\theta)$ of a tame elliptic maximal torus and a non\-singular character. Let $s \in S(F) \cap G(F)_\tx{reg}$ be topologically semi\-simple modulo center. The value of $\Theta_{\pi_{(S,\theta)}}$ at $s$ is given by
\[ e(G)\epsilon_L(X^*(T_G)_\C-X^*(S)_\C,\Lambda)D(s)^{-\frac{1}{2}}\sum_{w \in N(S,G)(F)/S(F)}\Delta_{II}^\tx{abs}[a,\chi'']({^w}s)\theta({^w}s). \]
\end{thm}
To briefly explain the notation, $D(s)=|\prod_\alpha (1-\alpha(s))|$ is the usual Weyl discriminant, the product being taken over all absolute roots of $S$, $e(G)$ is the Kottwitz sign of $G$ as in \cite{Kot83}, $T_G$ is the minimal Levi subgroup of the quasi\-split inner form of $G$, $\Lambda$ is an arbitrarily chosen non\-trivial character of the additive group of the base field $F$, $\epsilon_L$ is the root number of the given virtual Artin representation of degree $0$, and $\Delta_{II}^\tx{abs}$ is the function of $S(F)$ given by the formula
\[ \Delta_{II}^\tx{abs}(s) = \prod_\alpha \chi''_\alpha\left(\frac{\alpha(s)-1}{a_\alpha}\right). \]
The product runs over the $\Gamma$-orbits of absolute roots of $S$ that are symmetric, i.e. invariant under multiplication by $-1$, and $\Gamma$ is the absolute Galois group of $F$. If $\alpha$ represents such an orbit, we can associate the subgroups $\Gamma_\alpha \subset \Gamma_{\pm\alpha} \subset \Gamma$ and the corresponding field extensions $F_\alpha/F_{\pm\alpha}/F$. Then $a_\alpha \in F_\alpha^\times$ and $\chi''_\alpha : F_\alpha^\times \to \C^\times$ are computed explicitly in terms of $\theta$, and $a_\alpha$ depends moreover on $\Lambda$. We refer to \cite[\S4]{FKS} for the precise formulas. The first main takeaway is this:

\begin{center}
\textit{All constituents of this formula make sense for $F=\R$, and with this interpretation this formula recovers Harish-Chandra's formula from Theorem \ref{thm:hc2}.} 
\end{center}
There is however a key difference: Theorem \ref{thm:k3} applies only to very special elements of $S(F) \cap G(F)_\tx{reg}$ -- those that are topologically semi\-simple modulo center. It may happen that there are no such elements at all! So the \emph{values} of the given function may not uniquely characterize the representation $\pi_{(S,\theta)}$. There is a sense in which the \emph{formula} itself does characterize it, but such a statement may be met with skepticism by some colleagues, and in any event the question remains as to characterizing $\pi_{(S,\theta)}$ by its character function. There are two approaches to this problem. One, taken by Chan--Oi in \cite{CO21}, is to extend the validity of this formula to some more general elements of $S(F) \cap G(F)_\tx{reg}$ and prove that the resulting values are enough to characterize the representation; so far this has been successful under additional assumptions on $(S,\theta)$, including the assumption that $S$ is unramified. One can hope that such methods can be generalized to yield the validity of Theorem \ref{thm:k3} for all elements of $S(F)$ whose topologically semi-simple modulo center part is regular. The other approach, taken by \cite[\S4]{FKS}, is to establish a more general character formula, valid for all elements of $G(F)_\tx{reg}$ and all $(S,\theta)$, but under stricter conditions on $F$, as follows.

\begin{thm}[{\cite[Theorem 4.3.5]{FKS}}] \label{thm:k4}
Assume $F$ has characteristic zero and $p$ does not divide the order of the Weyl group of $G$ and is larger than $(2+e)n$, where $e$ is the ramification degree of $F/\Q_p$ and $n$ is the smallest dimension of a faithful algebraic representation of $G$. For any $\gamma \in G(F)_\tx{reg}$ with topological Jordan decomposition modulo center $\gamma = \gamma_0 \cdot \gamma_{0+}$ the value of $D(\gamma)^{\frac{1}{2}}\Theta_{\pi_{(S,\theta)}}(\gamma)$ is given by
\[ 
e(G)e(J)\epsilon_L(X^*(T_{G})_\C-X^*(T_{J})_\C,\Lambda)\!\!\!\!\!\!\!\!\!\!\!\!\!\!\!\!
\sum_{\substack
	{g \in S(F) \lmod G(F)/J(F) \\
	{^g}{\gamma_0} \in S(F)}
}\!\!\!\!\!\!\!\!\!\!\!\!\!\!\!\!
\Delta_{II}^\tx{abs}[a,\chi'']({^g}{\gamma_0}) \theta({^g}{\gamma_0})
\widehat O_{X^{g}}^{J}(\log \gamma_{0+}).
\]
\end{thm}
To explain the new notation, $J$ is the identity component of the centralizer of $\gamma_0$ in $G$, $X$ is any element of $\tx{Lie}^*(S)(F)$ such that $\theta(\exp(Y))=\Lambda(\<X,Y\>)$ for all $Y \in \tx{Lie}(S)(F)_{0+}$, and $\hat O_{X^g}^J$ is the (renormalized) function on $\tx{Lie}(J)(F)$ representing the Fourier transform of the orbital integral on $\tx{Lie}^*(J)(F)$ at $X^g$. In $\Delta_{II}^\tx{abs}$ we now drop those roots trivial on $^g\gamma_0$.

The reason we impose the stricter conditions on $F$ is so that the exponential map converges on the set of topologically nilpotent elements of $\tx{Lie}(G)(F)$, cf. \cite[App. A]{DR09}. This leads naturally to the following question:

\begin{qst} \label{qst:green}
Is there a formulation of the above character formula in which the function
$\widehat{O}_{X^g}^J(\log \gamma_{0+})$ is replaced by another function
on the set of topologically unipotent modulo center elements of
$J(F)$, which does not involve the logarithm map, but is still conjugation-invariant.
Such a function can be seen as a $p$-adic analog of Lusztig's Green functions, and this formulation would be valid also in positive characteristic and with possibly weaker conditions on $p$, as it avoids the use of the logarithm.
\end{qst}

It should be noted that the results of \cite{Spice21} are formulated with weaker conditions on $F$, but use a pseudo\-logarithm map that may not have good equivariance properties.

\begin{rem}
For Theorems \ref{thm:k3} and \ref{thm:k4} one must use the twisted Yu construction of \cite{FKS}, which is a modification of the original Yu construction whose purpose is to remedy an error in \cite{Yu01} that goes back to \cite{gerardin:weil}. That error invalidates some results of \cite{Yu01}, rendering invalid Yu's proof that the construction produces irreducible supercuspidal representations. It was shown by Fintzen in \cite{Fin19} that despite the error, the original Yu construction does produce irreducible supercuspidal representations. Nonetheless, the error introduces problems that lead to the appearance of auxiliary sign characters in the character formula (the characters $\epsilon^\tx{ram}$ and $\epsilon_{f,\tx{ram}}$ of \cite[Corollaries 4.8.2, 4.10.1]{KalRSP}, as well as the character of \cite[Proposition 5.27]{KalDC}), which make some applications of the resulting formula very difficult. In addition, the arguments of \cite{Spice21} can only be carried out for the twisted Yu construction.
\end{rem}

\subsection{Double covers of tori} \label{sub:dc}

In the archimedean setting Adams and Vogan \cite{AV92} have shown that Theorem \ref{thm:hc1}, as well as the local Langlands correspondence that will be our next topic, are more naturally formulated if instead of characters $\theta$ of an elliptic maximal torus $S(\R)$ one uses genuine characters of a certain double cover $S(\R)_\rho$. This double cover is obtained by choosing a Borel $\C$-subgroup $B$ of $G$ that contains $S$ and considering the algebraic double cover $S_\rho$ obtained as the pull\-back of the diagram $S \stackrel{2\rho}{\lrw} \C^\times \stackrel{2}{\llw} \C^\times$, where $2\rho$ is the sum of the $B$-positive absolute roots and $2$ denotes the squaring map. Then $S(\R)_\rho$ is defined as the preimage of $S(\R)$ under the isogeny $S_\rho(\C) \to S(\C)$. Note that there is a canonical character $\rho : S_\rho \to \mb{G}_m$. The choice of $B$ is immaterial, as one can take the limit over all possible choices. 

One can then reformulate Theorems \ref{thm:hc1} and \ref{thm:hc2} as the statement that there is a bijection between the set of discrete series representations of $G(\R)$ and the set of $G(\R)$-conjugacy classes of pairs $(S,\tilde\theta)$, where $S \subset G$ is an elliptic maximal torus and $\tilde\theta$ is a genuine character of the double cover $S(\R)_\rho$ such that $d\tilde\theta$ is regular. Note that we are not restricting here to what we called ``regular'' discrete series representations, i.e. this formulation of the theorem covers all discrete series representations. Moreover, the representation $\pi_{(S,\tilde\theta)}$ is the unique one whose Harish-Chandra character function evaluated at an element $s \in S(\R) \cap G(\R)_\tx{reg}$ has the form
\begin{equation} \label{eq:chardc-r}
(-1)^{q(G)}\frac{\sum\limits_{w \in N(S,G)(\R)/S(\R)}\tx{sgn}(w)\tilde\theta(\tilde s^w)}{\prod_{\alpha>0}(\alpha^\frac{1}{2}(\tilde s)-\alpha^{-\frac{1}{2}}(\tilde s))},
\end{equation}
where $\alpha>0$ runs over all absolute roots that are positive with respect to the Weyl chamber determined by the regular element $d\tilde\theta$, and $\tilde s \in S(\R)_\rho$ is any lift of $s$. Both numerator and denominator are well-defined genuine functions on $S(\R)_\rho$, and their quotient descends to $S(\R)$.

The double cover of Adams--Vogan generalizes to all local fields, but the generalization takes a different form than the original definition, in that it is of Galois\-theoretic rather than algebraic nature. Without going into technical details, for which we refer to \cite{KalDC}, we just mention that for any local field $F$, a connected reductive $F$-group $G$, and a maximal torus $S \subset G$, there exists a double cover $S(F)_\pm$ whose elements can be represented by tuples $(s,(\delta_\alpha))$ with $s \in S(F)$ and $\delta_\alpha \in F_\alpha^\times$ for every symmetric $\alpha \in R(S,G)$ such that $\delta_{\sigma(\alpha)}=\sigma(\delta_\alpha)$ for all $\sigma \in \Gamma$ and $\delta_\alpha/\delta_{-\alpha}=\alpha(s)$. When $F=\R$ and the torus $S$ is elliptic then $S(\R)_\pm$ is canonically identified with $S(\R)_\rho$.

Theorems \ref{thm:k1} and \ref{thm:k3} take the following shape in terms of this double cover: There is a natural bijection between the set of $G(F)$-conjugacy classes of regular supercuspidal representations and the set of $G(F)$-conjugacy classes of pairs $(S,\tilde\theta)$, where $S$ is a tame elliptic maximal torus and $\tilde\theta$ is a regular genuine character of the double cover $S(F)_\pm$. For any $s \in S(F) \cap G(F)_\tx{reg}$ that is topologically semi\-simple modulo center $\Theta_{\pi_{(S,\tilde\theta)}}$ takes the value
\[ e(G)\epsilon_L(X^*(T_G)_\C-X^*(S)_\C,\Lambda)D(s)^{-\frac{1}{2}} \sum_{w \in N_G(S)(F)/S(F)} a_S(\tilde s^w)\tilde\theta(\tilde s^w), \]
where $a_S : S(F)_\pm \to \{\pm 1\}$ is the genuine function sending $\tilde s=(s,(\delta_\alpha)) \in S(F)_\pm$ to
\[ \prod_\alpha \kappa_\alpha\left(\frac{\delta_\alpha-\delta_{-\alpha}}{a_\alpha}\right), \]
the product runs over the set of $\Gamma$-orbits of symmetric elements in $R(S,G)$, and $\kappa_\alpha : F_{\pm\alpha}^\times \to \{\pm 1\}$ is the quadratic character associated to the extension $F_\alpha/F_{\pm\alpha}$. One can also formulate Theorem \ref{thm:k4} in terms of $S(F)_\pm$; we skip this for now, but will formulate an analogous formula when discussing the local Langlands correspondence.

The advantage of using $S(F)_\pm$ is that it removes the somewhat mysterious characters $\chi''_\alpha$ used in Theorem \ref{thm:k3} (or rather it clarifies their role as mediating between characters of $S(F)$ and genuine characters of $S(F)_\pm$). Unfortunately, unlike in the archimedean case, this formulation does not allow the parameterization in terms of pairs $(S,\tilde\theta)$ to be extended beyond the case of regular supercuspidal representations.

\section{The local Langlands correspondence}

\subsection{The basic version}

Let $F$ be a local field and $G$ be a connected reductive $F$-group. Let $^LG=\hat G \rtimes \Gamma$ be the $L$-group of $G$ and let $L_F$ be the local Langlands group of $F$, i.e. the Weil group $W_F$ when $F$ is archimedean, or the group $W_F \times \tx{SL}_2(\C)$ when $F$ is non\-archimedean. The basic version of the local Langlands conjecture states that there exists a surjective finite\-to\-one map from the set of equivalence classes of irreducible admissible representations of $G(F)$ to the set of $\hat G$-conjugacy classes of relevant $L$-parameters  $\varphi : L_F \to {^LG}$. The fiber over $\varphi$ is called an \emph{$L$-packet}, denoted by $\Pi_\varphi(G)$.

There are reduction steps on the side of $L$-parameters that are parallel to the reduction steps ``admissible'' $\longrightarrow$ ``tempered'' $\longrightarrow$ ``discrete'' in the classification of irreducible admissible representations, but amount to simple exercises. The step ``admissible'' $\longrightarrow$ ``tempered'' produces from an arbitrary Langlands parameter $\varphi$ a triple $(P,\varphi_M,\nu)$ consisting of a parabolic subgroup $P$ of $G$, a tempered parameter $\varphi_M$ for the Levi quotient $M$ of $P$, and an element $\nu$ of $X^*(A_M)_\R$ that lies in the $P$-positive open cone, cf. \cite{SZ18} for the non\-archimedean case. The $L$-packet $\Pi_\varphi(G)$ then consists of the representations $j_P^G(\sigma\otimes\nu)$ for any $\sigma \in \Pi_{\varphi_M}(M)$. The step ``tempered'' $\longrightarrow$ ``discrete'' is even simpler, and just records the Levi subgroup $M$ of $G$ so that a given tempered $L$-parameter factors through $^LM$ and through no smaller Levi subgroup. The $L$-packet $\Pi_\varphi(G)$ consists of the irreducible constituents of the parabolic induction $i_P^G(\sigma)$ for any $\sigma \in \Pi_{\varphi}(M)$.

This reduces the construction of the correspondence to the case of discrete parameters, i.e. those that do not factor through any proper Levi subgroup, and essentially discrete series representations. At this point it becomes clear that the conjecture, as stated so far, is almost vacuous: nothing prevents us from randomly assigning discrete series representations to discrete parameters. This raises the following fundamental question, raised on various occasions by M. Harris, K. Buzzard, and others, which is so far unresolved in full generality:

\begin{qst} Find a list of properties that uniquely characterize the local Langlands correspondence.
\end{qst}
As discussed above, it is enough to answer this question for discrete parameters. While eventually a compatibility with a given global correspondence would be a key requirement, at the moment this is not feasible, and we seek a purely local characterization.

A number of expected properties have already been formulated, for example compatibility with central and cocentral characters and homomorphisms with abelian kernel and cokernel (\cite[10.3(1),(2),(5)]{BorCor}), the strong tempered $L$-packet conjecture (a strengthening of \cite[Conjecture 9.4]{Sha90} stating that each tempered $L$-packet on a quasi-split group contains a unique member that is generic with respect to a fixed Whittaker datum), the formal degree conjecture (\cite{HII08}), and the contragredient conjecture (\cite{AV16}, \cite{KalGen}). These are however not enough to pin down the correspondence uniquely. The following property is also expected:

\begin{cnj} \label{cnj:as}
Each discrete series $L$-packet is atomically stable, i.e. there exists a linear combination of the Harish-Chandra characters of its members that is a stable distribution, and no proper subset of the $L$-packet has this property.
\end{cnj}
\footnotetext{Conjecture \ref{cnj:as} is expected more generally for tempered $L$-packets, but not for non-tempered $L$-packets; the latter need to be enlarged to Arthur packets, or more generally ABV packets, in order to provide stable distributions, cf. \cite{ArtUARC}, \cite{ABV92}.}

It is expected that Conjecture \ref{cnj:as} uniquely characterizes the partition of the set of equivalence classes of irreducible discrete series representations of $G(F)$ into $L$-packets. However, it does not determine the matching between $L$-packets and $L$-parameters. 

In the case of $\tx{GL}_N$, stability is a vacuous condition and $L$-packets are singletons. On the other hand, Henniart has found \cite{Hen93,Hen02} a list of conditions that uniquely determine the local Langlands correspondence for $\tx{GL}_N$ when $F$ is non-archimedean. Besides the already listed conditions regarding central and cocentral characters and contragredient, what is needed is equality of $L$- and $\epsilon$-factors of pairs, which on the Galois side are the Artin factors of the tensor product of the two Galois representations, and on the automorphic side are given by Rankin--Selberg integrals. While analogous factors can be defined for some other groups as well, such as classical groups, it is unfortunately not yet known how to define them for general reductive groups intrinsically. For some interesting ideas in this direction, see \cite{BSN16}. Another approach to characterizing the local Langlands correspondence for non-archimedean $F$ was pioneered by Scholze in \cite{Scholze13} for the group $\tx{GL}_N$, and extended to a certain list of other groups in \cite{BMY20} based on \cite{SS13}. 

One way to characterize the assignment $\varphi \mapsto \Pi_\varphi(G)$ would be to associate to each $\varphi$ a stably invariant distribution $S\Theta_\varphi^G$ that would be the stable character of the corresponding $L$-packet (unique up to non\-zero scalar multiple by Conjecture \ref{cnj:as}). In the archimedean case, where the local correspondence has been constructed by Langlands \cite{Lan89}, this stable distribution (in fact function) can be described most conceptually using the double covers discussed in \S\ref{sub:dc}. According to Adams--Vogan \cite{AV92}, this description is as follows. The existence of a discrete parameter $\varphi$ easily implies the existence of an elliptic maximal torus $S \subset G$. There is an $L$-group $^LS_\pm$ associated to the double cover $S(\R)_\pm$. One key property of the double cover is that there is a canonical $\hat G$-conjugacy class of $L$-embeddings $^LS_\pm \to {^LG}$ (this is not the case for the $L$-group $^LS$ of the torus $S$ itself), cf. \cite[\S4.1]{KalDC}. It is again easy to see that $\varphi$ factors through this $L$-embedding, and thus leads canonically to a genuine character $\tilde\theta$ of $S(\R)_\pm$, well-defined up to $W_G(S)(\R)$. The stable character associated to $\varphi$ is uniquely characterized by its restriction to $S(\R)$, where it takes the form
\begin{equation} \label{eq:chardc-rst}
S\Theta_\varphi^G(s)=
(-1)^{q(G^*)}\cdot \frac{\sum\limits_{w \in W_G(S)(\R)}\tx{sgn}(w)\tilde\theta(\tilde s^w)}{\prod_{\alpha>0}(\alpha^\frac{1}{2}(\tilde s)-\alpha^{-\frac{1}{2}}(\tilde s))}.
\end{equation}
Here again $\alpha>0$ means $\<\alpha^\vee,d\tilde\theta\>>0$. We denote by $W_G(S)=N_G(S)/S$ the absolute Weyl group, and by $G^*$ the quasi-split inner form of $G$. Note the very close relationship to \eqref{eq:chardc-r}. It may appear odd that we insist on the constant $(-1)^{q(G^*)}$ even though the entire function is supposed to be well-defined only up to a constant; we will see in the next subsection that in fact this function is conjecturally well-defined ``on the nose'', and not just up to a constant.

Turning to a non-archimedean base field $F$, we can use these ideas of Adams--Vogan and our experience from Theorems \ref{thm:k3} and \ref{thm:k4} to formulate the following conjectures describing the stable character associated to a \emph{supercuspidal} parameter $\varphi$, i.e. a discrete parameter that is trivial on the subgroup $\tx{SL}_2(\C)$ of $L_F$. It is easy to see that, when $G$ is tame and $p$ does not divide the order of the Weyl group of $G$, such a parameter determines a stable conjugacy class of elliptic maximal tori $S \subset G$, and that $\varphi$ factors through the canonical embedding $^LS_\pm \to {^LG}$, thereby providing a genuine character $\tilde\theta$ of $S(F)_\pm$.

\begin{cnj} \label{cnj:schar1}
Let $\gamma \in G(F)$ be regular semi-simple and topologically semi-simple modulo center. Then $S\Theta_\varphi^G(\gamma)$ is zero unless $\gamma$ lies in the image of an admissible embedding $S \to G$, in which case (after identifying $S$ with that image), we have	
\begin{equation} \label{eq:schar1}
S\Theta_\varphi^G(\gamma) = \epsilon_L(X^*(T_G)_\C-X^*(S)_\C)D(\gamma)^{-\frac{1}{2}}\sum_{w \in W_G(S)(F)}[a_S \cdot \tilde\theta](\gamma^w).
\end{equation}
\end{cnj}

Note the strong similarity between \eqref{eq:chardc-rst} and \eqref{eq:schar1}. In fact, this is more than just a similarity: 

\begin{center}
\textit{Formula \eqref{eq:schar1} makes sense for any local field $F$, and recovers Formula \eqref{eq:chardc-rst} when $F=\R$. Therefore, it gives a conjectural description of the stable character associated to a discrete Langlands parameter $\varphi : W_F \to {^LG}$ that is uniform for any local field $F$.
}
\end{center}
Of course, this conjecture only applies to discrete parameters that, in the non-archimedean case, have trivial restriction to $\tx{SL}_2(\C)$, and in addition $p$ is prime to $|W_G(S)|$. This conjecture was the guiding principle behind the constructions of \cite{KalRSP} and \cite{KalSLP}. One drawback it has is that in the non-archimedean case there may not be enough topologically semi-simple elements of $S(F)$ to fully determine the function $S\Theta_\varphi^G$. This is not an issue in the setting of loc. cit., because we are not using the values of the function, but rather the \emph{entire formula}, which carries more information. Nonetheless, a more complete solution is desirable. It is conceivable that the ideas of \cite{CO21} might lead to a stronger version of this formula. Another approach is to allow arbitrary regular semi-simple elements of $G(F)$, in the vein of Theorem \ref{thm:k4}. This leads to the following conjecture.

\begin{cnj}[{\cite[\S4.4]{KalDC}}] \label{cnj:schar}
For any strongly regular semi\-simple $\gamma \in G(F)$ with topological Jordan decomposition modulo center $\gamma=\gamma_0 \cdot \gamma_{0+}$,
\[ S\Theta_\varphi^G(\gamma) = e(J)\epsilon_L(X^*(T_G)_\C-X^*(T_J)_\C)D(\gamma)^{-\frac{1}{2}}\sum_{j : S \to J}[a_S \cdot \tilde\theta](\gamma_0^j) \cdot \widehat{SO}^J_{^jX}(\log(\gamma_{0+})), \]
assuming $F$ has characterstic zero and $p \geq (2+e)n$. 
\end{cnj}

The notation is the same as that in Theorem \ref{thm:k4}, except now we are using the stable orbital integral at $^jX$ instead of the usual orbital integral, and the sum runs over the set of stable classes of admissible embeddings $S \to J$.  

Again, the reason we require the characteristic of $F$ to be zero and $p$ to be very large is to ensure the convergence of $\exp$ on $\tx{Lie}(G)(F)_{0+}$, in particular on $\tx{Lie}(J)(F)_{0+}$. A positive resolution to Question \ref{qst:green} would weaken this requirement.

The most general constructions of the non-archimedean basic local Langlands correspondence are given by Genestier--Lafforgue \cite{GL17} in positive characteristic and Fargues--Scholze \cite{FS21} for arbitrary non\-archimedean local fields. These constructions only produce semi\-simplified parameters, but, at least in positive characteristic, recent work of Gan--Harris--Sawin \cite{GHS} (based on arguments of Gan--Lomel\'i \cite{GanLomeli18} in the case of classical groups) provides a unique enrichment of such a semi\-simplified parameter to a full Langlands parameter when the representation in question is supercuspidal. Earlier constructions include \cite{LRS93,HT01,Hen00} for $\tx{GL}_N$ and \cite{Art13} for quasi-split symplectic and orthogonal groups in characteristic zero, and \cite{HKT19} for generic supercuspidal representations of the exceptional group $G_2$.

At the moment there are many open questions regarding the constructions \cite{GL17} and \cite{FS21}, such as
\begin{enumerate}%
	\item Is the map $\pi \mapsto \varphi$ surjective?
	\item Are the resulting $L$-packets always finite?
	\item Does the construction of Fargues--Scholze specialize to that of Genestier--Lafforgue when $F$ has positive characteristic?
	\item Do conjectures \ref{cnj:as} and \ref{cnj:schar} hold?
	\item Does the formal degree conjecture hold?
\end{enumerate}

On the other hand, \cite{KalSLP} gives an \emph{explicit} construction of the correspondence under the assumption that $G$ is tame and $p$ does not divide the order of the Weyl group, and $\varphi$ is a supercuspidal parameter, building on prior work \cite{DR09}, \cite{Ree08}, \cite{KalEpi}. This setting is more restrictive than that of \cite{FS21} or \cite{GL17}, but in return provides much more knowledge about the resulting correspondence. For example, we know that
\begin{enumerate}%
	\item The map $\pi \mapsto \varphi$ has as domain all non\-singular supercuspidal representations, and as image all supercuspidal parameters.
	\item The map $\pi \mapsto \varphi$ is compatible with central and cocentral characters.
	\item The resulting $L$-packets are always finite, and in fact have the desired internal structure (see next section).
	\item Both conjectures \ref{cnj:as} and \ref{cnj:schar} hold (\cite[\S4.4]{FKS}).
	\item The formal degree conjecture holds, as shown by Schwein \cite{Schwein21} and Ohara \cite{Ohara21}.
\end{enumerate}

The question whether the constructions of \cite{FS21} and \cite{KalSLP} agree is equivalent (in the setting of $F$ having characteristic zero and $p$ being sufficiently large) to the question of whether \cite{FS21} satisfies Conjecture \ref{cnj:schar}. A strong indication that they agree would be given if Conjecture \ref{cnj:schar1} could be proved instead; the latter can also be pursued for \cite{GL17}, since it allows the characteristic of $F$ to be positive.

The opposite setting of that of supercuspidal parameters and non-singular (i.e. semi-simple) supercuspidal representations is that of unipotent supercuspidal representations, and more generally arbitrary unipotent representations. Much progress has been made on the local Langlands correspondence for these representations via detailed study of affine Hecke algebras and formal degrees \cite{Lus95,Lus02}, \cite{Ree00}, \cite{FOS20}, \cite{Soll19}.

\subsection{The refined version}

For many applications, such as the Gan--Gross--Prasad conjecture, or the multiplicity formula for discrete automorphic representations, the basic version of the local Langlands conjecture is insufficient, because it describes packets of representations rather than individual representations. The refined local Langlands conjecture remedies this by enhancing the notion of a Langlands parameter to allow the description of individual irreducible admissible representations. In fact, already the statement of the Hiraga--Ichino--Ikeda conjecture requires the refined correspondence, a point that we glossed over in the previous subsection. 

As pointed out by Vogan \cite{Vog93}, the refined conjecture requires a rigidification of the concept of inner forms of reductive groups. In the archimedean case, a good rigidification was obtained by Adams--Barbasch--Vogan \cite{ABV92}. In the non-archimedean case, different ways of rigidification have led to different versions of the refined conjecture. We only give a brief summary, referring the reader to \cite{KalSimons} for more details. 

The set of equivalence classes of inner forms of a connected reductive group $G$ is given by $H^1(F,G_\tx{ad})$, where $G_\tx{ad}=G/Z(G)$ is the adjoint group of $G$. In every inner class there is a unique quasi-split form, and we normalize things by taking $G$ to be that form. A rigidification of the notion of an inner form can be achieved by choosing a Galois gerbe $\mc{E}$. When $F$ has characteristic zero such a gerbe can be understood, following Langlands--Rapoport \cite{LR87}, as an extension $1 \to u(\ol{F}) \to \mc{E} \to \Gamma \to 1$ of the absolute Galois group $\Gamma$ of $F$ by an algebraic or pro-algebraic group $u$. Following Kottwitz \cite{Kot97}, one then considers the set $H^1_\tx{bas}(\mc{E},G)$ of cohomology classes of $\mc{E}$ with values in $G(\ol{F})$ whose restriction to $u$ factors through the center $Z(G)$ of $G$. There is a natural embedding $H^1(F,G) \to H^1_\tx{bas}(\mc{E},G)$ and a natural map 
\begin{equation} \label{eq:gerbe}
H^1_\tx{bas}(\mc{E},G) \to H^1(F,G_\tx{ad}).	
\end{equation}
A rigidification of an inner form of $G$ is the choice of an element of $H^1_\tx{bas}(\mc{E},G)$ that lifts the class of that inner form.

Vogan's notion of pure inner forms comes from the trivial Galois gerbe $\mc{E}^\tx{triv}=\Gamma$, for which $u=\{1\}$. Kottwitz's theory of isocrystals with $G$-structure \cite{Kot85} employs the gerbe $\mc{E}^\tx{iso}$ of the Tannakian category of $F$-isocrystals. The problem with $\mc{E}^\tx{triv}$ and $\mc{E}^\tx{iso}$ is that in general the map \eqref{eq:gerbe} is not surjective, so not all inner forms can be rigidified. For $\mc{E}^\tx{iso}$ that map is surjective when $Z(G)$ is connected. In \cite{KalRI} I define a gerbe $\mc{E}^\tx{rig}$ for which the map \eqref{eq:gerbe} is always surjective. It turns out that, when $F=\R$, this gerbe recovers the notion of strong real forms introduced by \cite{ABV92}. A key property of the gerbes $\mc{E}^\tx{iso}$ and $\mc{E}^\tx{rig}$ is that they satisfy a generalization of Tate--Nakayama duality.

When $F$ has positive characteristic the simplified concept of a Galois gerbe as an extension of the absolute Galois group becomes inadequate, due to the possible non-smoothness of $u$. Despite this difficulty, Peter Dillery \cite{Dillery20} has found a way to construct a suitable analog of $\mc{E}^\tx{rig}$. In fact, his construction works uniformly for all non-archimedean local fields and recovers $\mc{E}^\tx{rig}$ when $F$ has characteristic zero. Therefore, we now have a satisfactory definition of $\mc{E}^\tx{rig}$ for any local field.

The refined local Langlands conjecture parameterizes all irreducible admissible representations of all inner forms of $G$ at once. More precisely, one considers tuples $(G',\xi,z,\pi)$, where $G'$ is a connected reductive $F$-group, $\xi : G_{F^s} \to G'_{F^s}$ is an isomorphism of $F^s$-groups, where $F^s$ is a fixed separable closure of $F$, $z \in Z^1_\tx{bas}(\mc{E}^\tx{rig},G)$, and $\xi^{-1}\sigma(\xi)=\tx{Ad}(\bar z_\sigma)$, where $\bar z_\sigma$ is the image of $z$ in $Z^1(F,G_\tx{ad})$ (since $G_\tx{ad}$ is smooth we can interpret the latter as \'etale, i.e. Galois, cohomology), and finally $\pi$ is an irreducible admissible representation of $G'(F)$. An isomorphism $(G_1,\xi_1,z_1,\pi_1) \to (G_2,\xi_2,z_2,\pi_2)$ is a pair $(g,f)$ with $g \in G(F^s)$ and $f : G_1 \to G_2$ an isomorphism of $F$-groups such that $f \circ \xi_1 = \xi_2 \circ \tx{Ad}(g)$ and $z_2(e)=gz_1(e)\sigma_e(g)^{-1}$; here $\sigma_e \in \Gamma$ is the image of $e \in \mc{E}^\tx{rig}$. The key property of the set $Z^1_\tx{bas}(\mc{E}^\tx{rig},G)$ is that if we fix the triple $(G',\xi,z)$, then two tuples $(G',\xi,z,\pi_1)$ and $(G',\xi,z,\pi_2)$ are isomorphic if and only if the representations $\pi_1$ and $\pi_2$ of $G'(F)$ are equivalent.

Assuming the validity of the basic local Langlands conjecture, we can define for each Langlands parameter $\varphi : L_F \to {^LG}$ the ``compound $L$-packet'' $\Pi_\varphi$ as the set of isomorphism classes of tuples $(G',\xi,z,\pi)$, for all possible $(G',\xi,z)$ and all $\pi \in \Pi_\varphi(G')$. We further let $S_\varphi$ be the centralizer of $\varphi$ in $\hat G$, and $S_\varphi^+$ and $Z([\hat{\bar G}]^+)$ be the preimages of $S_\varphi$ and $Z(\hat G)^\Gamma$ in the universal cover $\hat{\bar G}$ of $\hat G$. The refined conjecture, inspired by the work of Adams-Barbasch-Vogan and Vogan, is then the following.

\begin{cnj}[{\cite[\S5.4]{KalRI}}] \label{cnj:refined}
Fix a Whittaker datum $\mf{w}$ for $G$.
\begin{enumerate}%
	\item There exists a bijection $\iota_{\varphi,\mf{w}}$ that is the top map in the following commutative diagram
	\[ 
	\xymatrix{
	\Pi_\varphi\ar[r]^-{\iota_{\varphi,\mf{w}}}\ar[d]&\tx{Irr}(\pi_0(S_\varphi^+))\ar[d]\\
	H^1_\tx{bas}(\mc{E}^\tx{rig},G)\ar[r]&\pi_0(Z(\hat{\bar G})^+)^*
	}
	\]
	in which the left map sends the isomorphism class of $(G',\xi,z,\pi)$ to the class of $z$, the right map is the central character map, and the bottom map is generalized Tate--Nakayama duality. This bijection relates the unique $\mf{w}$-generic constituent of $\Pi_\varphi$ to the trivial representation of $\pi_0(S_\varphi^+)$.
\end{enumerate}
Given a semi-simple element $\dot s \in S_\varphi^+$, a rigid inner twist $(G',\xi,z)$, and a \emph{tempered} parameter $\varphi$, define the virtual character
\[ \Theta^{G',\xi,z}_{\varphi,\mf{w},\dot s} = e(G')\sum_{\pi \in \Pi_\varphi(G')} \tx{tr}(\iota_{\varphi,\mf{w}}(G',\xi,z,\pi)(\dot s))\cdot \Theta_\pi. \]
\begin{enumerate}%
	\item The distribution $\Theta^{G',\xi,z}_{\varphi,\mf{w},1}$ is stable and independent of $\mf{w}$ and $z$.
	\item For a general $\dot s \in S_\varphi^+$, the distribution $\Theta^{G',\xi,z}_{\varphi,\mf{w},\dot s}$ is the endoscopic lift of the distribution $\Theta^{H,\tx{id},1}_{\varphi,*,\dot s}$ for the endoscopic datum $(H,\dot s)$ associated to $\varphi$, with respect to the transfer factor normalized via $\mf{w}$ and $z$ as in \cite[(5.10)]{KalRI}.
\end{enumerate}
\end{cnj}

\begin{rem} \label{rem:refined}
\ \\[-20pt]
\begin{enumerate}%
	\item Point (3) specifies the bijection $\iota_{\varphi,\mf{w}}$ uniquely, provided such a bijection exists. 
	\item The distribution in (2) is what we referred to as $S\Theta_\varphi^{G'}$ in the previous subsection. Note that the inner twist $\xi$ is used to identify $^LG'$ with $^LG$, so if we use $^LG$ as codomain for Langlands parameters, we should write $S\Theta_\varphi^{G',\xi}$ to indicate that this distribution depends also on $\xi$, not just $G'$.
	\item The fiber over $[z] \in H^1_\tx{bas}(\mc{E}^\tx{rig},G)$ of the left vertical map in (1) is by definition the $L$-packet $\Pi_\varphi(G_{\bar z})$ of the inner form $G_{\bar z}$ of $G$ associated to the image $[\bar z] \in H^1(F,G_\tx{ad})$ of $z$. Note that there can be two distinct $[z_1],[z_2] \in H^1(\mc{E}^\tx{rig},G)$ mapping to $[\bar z]$. This leads to the appearance of the same $L$-packet $\Pi_\varphi(G_{\bar z})$ multiple times in the compound $L$-packet $\Pi_\varphi$. This ``overcounting'' is the spectral incarnation of the rigidification of the notion of inner forms. 

	\item When $F$ is non-archimedean, the bottom map in (1) is bijective. This means that the $L$-packet $\Pi_\varphi(G_{\bar z})$ for the individual group $G_{\bar z}$ is described precisely by the set $\tx{Irr}(\pi_0(S_\varphi^+),[z])$ of irreducible representations of $\pi_0(S_\varphi^+)$ that transform under the action of $\pi_0(Z(\hat{\bar G})^+)$ via the character corresponding to $[z]$. As just discussed, there can be multiple $[z]$ lifting $[\bar z]$, leading to multiple ways to parameterize $\Pi_\varphi(G_{\bar z})$. Thankfully, it is a rather straightforward matter to relate these two parameterizations of $\Pi_\varphi(G_{\bar z})$, as we will discuss in the next subsection.

	\item When $F=\R$, the bottom map in (1) need not be injective nor surjective. These failures are well understood, cf \cite[\S3.4, \S4, Proposition 5.3]{KalRI}. The non-surjectivity implies that, for some $\eta \in \pi_0(Z(\hat{\bar G})^+)^*$, the set $\tx{Irr}(\pi_0(S_\varphi^+),\eta)$ does not parameterize any $L$-packet. The non-injectivity implies that for some $\eta \in \pi_0(Z(\hat{\bar G})^+)^*$, the set $\tx{Irr}(\pi_0(S_\varphi^+),\eta)$ parameterizes a union of $L$-packets over certain rigid inner forms. The set of rigid inner forms appearing in such a union consists of exactly those inner forms of $G$ that are related to each other by Galois 1-cocycles valued in the simply connected covers of their adjoint groups; this set is termed a \emph{$K$-group} by Arthur in \cite[\S1]{Art99}.

	\item An analogous conjecture can be stated with $\mc{E}^\tx{iso}$ or $\mc{E}^\tx{triv}$ in place of $\mc{E}^\tx{rig}$. One has to then replace $\pi_0(Z(\hat{\bar G})^+)^*$ by $X^*(Z(\hat G)^\Gamma)$ or $\pi_0(Z(\hat G)^\Gamma)^*$, respectively, and $\pi_0(S_\varphi^+)$ by $S_\varphi^\natural = S_\varphi/(S_\varphi \cap \hat G_\tx{der})^\circ$ or $\pi_0(S_\varphi)$, respectively. There is a commutative diagram 
	\[ 
	\xymatrix{
	\pi_0(S_\varphi)&\ar[l]\pi_0(S_\varphi^+)\ar[r]&S_\varphi^\natural\\
	\pi_0(Z(\hat G)^\Gamma)\ar[u]&\ar[l]\pi_0(Z(\hat{\bar G})^+)\ar[u]\ar[r]&Z(\hat G)^\Gamma\ar[u]
	}
	\]
	that relates the three versions of the conjecture. Given characters $\xi_1 \in \pi_0(S_\varphi)^*$ and $\xi_2 \in X^*(Z(\hat G)^\Gamma)$ with common image $\xi \in \pi_0(Z(\hat{\bar G})^+)^*$, this diagram induces bijections $\tx{Irr}(\pi_0(S_\varphi),\xi_1) = \tx{Irr}(\pi_0(S_\varphi^+),\xi) = \tx{Irr}(S_\varphi^\natural,\xi_2)$. Therefore, the three versions of the conjecture differ only in the amount of inner forms of $G$ that they can reach, and the way in which these inner forms are overcounted by rigidification. For a more thorough comparison, cf. \cite[\S4.2]{KalRIBG}.
\end{enumerate}
\end{rem}

When $F=\R$, Conjecture \ref{cnj:refined} was established by Shelstad in a series of papers, especially \cite{She82,SheTE1,SheTE2,SheTE3}, cf. also \cite[\S5.6]{KalRI}. A geometric approach to Conjecture \ref{cnj:refined} was developed in \cite{ABV92}, which also produces and parameterizes Arthur packets, not just $L$-packets.

For non-archimedean local fields of characteristic zero, the restriction of Conjecture \ref{cnj:refined} to quasi-split classical groups was proved by Arthur \cite{Art13}.

For more general groups over non-archimedean local fields, we have the following.

\begin{thm}[{\cite[\S4.4]{FKS}}] \label{thm:refined}
Under the assumptions of Theorem \ref{thm:k4}, Conjecture \ref{cnj:refined} holds for the construction of regular supercuspidal $L$-packets of \cite{KalRSP}, and more generally for all supercuspidal $L$-packets constructed \cite{KalSLP}, but in that larger generality point (3) is proved only for certain elements $\dot s$ (the remaining elements $\dot s$ are work in progress).
\end{thm}
It is at the moment not known if the constructions of \cite{GL17} or \cite{FS21} satisfy Conjecture \ref{cnj:refined}. This would follow from Theorem \ref{thm:refined} once Conjecture \ref{cnj:schar} is established for these constructions, via comparison with \cite{KalSLP}, at least under the assumptions under which that theorem and conjecture are formulated. We note that \cite{FS21} works with the concept of inner forms rigidified via $\mc{E}^\tx{iso}$, because isocrystals with $G$-structure appear naturally in the theory of the Fargues--Fontaine curve, but the relationship between the two versions of the refined correspondence is well understood as per Remark \ref{rem:refined}(6).

The refined local Langlands conjecture can be combined with a global version of $\mc{E}^\tx{rig}$, developed in \cite{KalGRI} over number fields, and in \cite{Dillery21} over global function fields, to obtain a precise formulation of the conjectural multiplicity formula for discrete automorphic representations originally due to Kottwitz \cite[(12.3)]{Kot84}, cf. \cite[\S4.5]{KalGRI}, \cite[\S5.4]{Dillery21}. Special cases of this formula have been proved by O. Ta\"ibi, cf. \cite{TaibMult}. One can also obtain a global multiplicity conjecture using the global version of $\mc{E}^\tx{iso}$ defined by Kottwitz in \cite{KotBG}, under the assumption that the global group $G$ has connected center and satisfies the Hasse principle. Under those conditions, the two global multiplicity formulas coming from $\mc{E}^\tx{iso}$ and $\mc{E}^\tx{rig}$ are equivalent, cf. \cite{KT18}.

Another setting in which information about the refined local Langlands conjecture, more precisely part (1) of Conjecture \ref{cnj:refined}, has been obtained is that of unipotent representations, cf. \cite{Lus95,Lus02}, \cite{Ree00}, \cite{FOS20}, \cite{Soll19}.
 
\subsection{Compatibility properties}

The basic version of the local Langlands correspondence is expected to satisfy many compatibility properties, cf. \cite[\S10]{BorCor}. The refined version of the local Langlands correspondence allows one to formulate more precise compatibility properties. Some of them turn out to be formal consequences of the refined conjecture, while others have to be proved independently. 

Among the simplest such properties are those regarding the dependence of the parameterizing bijection $\iota_{\varphi,\mf{w}}$ on the Whittaker datum $\mf{w}$ and on the element of $H^1_\tx{bas}(\mc{E}^\tx{rig},G)$ lifting the class of $H^1(F,G_\tx{ad})$ that describes the relevant inner form.

\subsubsection{Whittaker data} 
If another Whittaker datum $\mf{w}'$ is chosen, there exists an element $g \in G_\tx{ad}(F)$ such that $\mf{w}' = g\mf{w}g^{-1}$. Denote by $(\mf{w}',\mf{w})$ the image of $g$ under the connecting homomorphism $H^0(F,G_\tx{ad}) \to H^1(F,Z(G_\tx{sc}))$. Local Tate duality identifies $H^1(F,Z(G_\tx{sc}))$ with the dual of $H^1(F,Z(\hat G_\tx{sc}))=H^1(L_F,Z(\hat G_\tx{sc}))$. A Langlands parameter $\varphi : L_F \to {^LG}$ induces an action of $L_F$ on $\hat G$ by conjugation via $\varphi$, which coincides with the usual action of $L_F$ on $Z(\hat G_\tx{sc})$. Via the connecting homomorphism $H^0(\varphi(L_F),\hat G_\tx{ad})) \to H^1(L_F,Z(\hat G_\tx{sc}))$ in the resulting long exact cohomology sequence we can pull back the character $(\mf{w}',\mf{w})$ to the group $S_\varphi/Z(\hat G)^\Gamma = S_\varphi^+/Z(\hat{\bar G})^+$. The following result is stated in \cite[Theorem 4.3]{KalGen} for real groups or quasi-split classical $p$-adic groups, but the proof actually shows that it is a formal consequence of Conjecture \ref{cnj:refined}:

\begin{thm} \label{thm:generic}
The validity of Conjecture \ref{cnj:refined} implies
\[ \iota_{\varphi,\mf{w}'}(\dot\pi) = \iota_{\varphi,\mf{w}}(\dot \pi) \otimes (\mf{w}',\mf{w})\quad\forall \dot\pi \in \Pi_\varphi. \]
\end{thm}

\subsubsection{Rigidifying data}
The question of changing the rigidifying element can be resolved in an analogous way. Consider an element $\bar z \in Z^1(F,G_\tx{ad})$ leading to the inner form $G_{\bar z}$. Let $z_1,z_2 \in Z^1_\tx{bas}(\mc{E}^\tx{rig},G)$ both lift $\bar z$. Then $z_2z_1^{-1} \in Z^1(\mc{E}^\tx{rig},Z(G))$ and we denote by $[z_2z_1^{-1}]$ its class. Given a Langlands parameter $\varphi : L_F \to {^LG}$ and $\pi \in \Pi_\varphi(G_{\bar z})$ we can consider the elements $\dot\pi_1=(G_{\bar z},\xi,z_1,\pi)$ and $\dot\pi_2=(G_{\bar z},\xi,z_2,\pi)$ of $\Pi_\varphi$, where $\xi$ denotes the identification of $G_{F^s}$ with $G_{\bar z,F^s}$. These elements both describe the representation $\pi$ of the group $G_{\bar z}(F)$, but are distinct elements of the compound packet $\Pi_\varphi$. This is the overcounting phenomenon mentioned in the previous subsection. To relate these two ``reflections'' of $\pi$, consider the exact sequence
\[ 1 \to \pi_1(\hat G) \to \hat{\bar G} \to \hat G \to 1 \]
equipped with the action of $L_F$ via conjugation by $\varphi$. It leads to the differential $d : S_\varphi^+ \to Z^1(F,\pi_1(\hat G))$ that factors through $\pi_0(S_\varphi^+)$. We denote by $-d$ the composition of this differential with the inversion automorphism of the abelian group $\pi_1(\hat G)$. It is shown in \cite[\S\S6.1,6.2]{KalRIBG} that local Tate duality generalizes to a duality between $H^1(\mc{E}^\tx{rig},Z(G))$ and $Z^1(F,\pi_1(\hat G))$. The element $[z_2z_1^{-1}]$ thus becomes a character of $Z^1(F,\pi_1(\hat G))$, which can be pulled back to $\pi_0(S_\varphi^+)$ by $-d$.

\begin{thm}[{\cite[\S6.3]{KalRIBG}}] \label{thm:central}
The validity of Conjecture \ref{cnj:refined} implies 
\[ \iota_{\varphi,\mf{w}}(\dot\pi_2) = \iota_{\varphi,\mf{w}}(\dot\pi_1) \otimes (-d)^*([z_2z_1^{-1}]).\]
\end{thm}

\begin{rem} \label{rem:rx}
The cohomological constructions in Theorem \ref{thm:generic} and \ref{thm:central} are very closely related: the composition of the differential $d : S_\varphi^+ \to Z^1(F,\pi_1(\hat G))$ with the natural projection $Z^1(F,\pi_1(\hat G)) \to H^1(F,\pi_1(\hat G))$ and the natural map $\pi_1(\hat G) \to \pi_1(\hat G_\tx{ad}) = Z(\hat G_\tx{sc})$ equals the composition of the natural projection $S_\varphi^+ \to H^0(\varphi(L_F),\hat G_\tx{ad})$ with the connecting homomorphism $H^0(\varphi(L_F),\hat G_\tx{ad}) \to H^1(L_F,Z(\hat G_\tx{sc}))=H^1(F,Z(\hat G_\tx{sc}))$. This is embodied in the commutative diagrams \cite[(6.1),(6.2),(6.6)]{KalRIBG}.
\end{rem}

\subsubsection{Contragredients}
We now discuss the compatibility of the local Langlands conjecture with respect to taking contragredients. To state it, we write $\hat C$ for the Chevalley involution of $\hat G$, well-defined up to conjugation, and by $^LC$ its extension to an automorphism of $^LG$. We write $\pi^\vee$ to denote the contragredient of a representation $\pi$, and given $\dot\pi=(G',\xi,z,\pi)$ we write $\dot\pi^\vee=(G',\xi,z,\pi^\vee)$.

\begin{cnj} \label{cnj:contra}
Assume Conjecture \ref{cnj:refined}. Let $\varphi : L_F \to {^LG}$ be a Langlands parameter with compound $L$-packet $\Pi_\varphi$. Then 
\begin{enumerate}%
	\item $\Pi_{^LC \circ\varphi} = \{\dot\pi^\vee\,|\, \dot\pi \in \Pi_\varphi\}$, 
	\item[(1')] $S\Theta_{^LC\circ\varphi}^{G',\xi}(\delta)=S\Theta_\varphi^{G',\xi}(\delta^{-1})$ whenever $\varphi$ is tempered,
	\item $\iota_{^LC \circ \varphi,\mf{w}^{-1}}(\dot\pi^\vee) = (\iota_{\varphi,\mf{w}}(\pi)\circ\hat C^{-1})^\vee$.
\end{enumerate}
\end{cnj}

Part (1) can be stated for $L$-packets on an individual group, and hence without assuming Conjecture \ref{cnj:refined}. In this form it was formulated by Adams--Vogan in \cite{AV16}, where it was also proved for $F=\R$. Wen-Wei Li proved that statement of (1) in \cite{WenWeiLi19} for the semi-simplified basic correspondence of \cite{GL17}. Part (1') is a variation of part (1) -- it implies part (1) via linear independence of characters, provided one assumes Conjecture \ref{cnj:as}. Over $F=\R$ parts (1) and (1') are in fact equivalent, since Conjectures \ref{cnj:as} and \ref{cnj:refined} are known and $S\Theta_\varphi^{G',\xi}$ is simply the sum of $\Theta_\pi$ over all $\pi \in \Pi_\varphi(G')$.

Part (2) was formulated by D. Prasad, cf. \cite[Conjecture 2]{Dipendra19}. It was proved in \cite[Theorem 5.9]{KalGen} that Conjecture \ref{cnj:refined} for $G$ together with part (1') for all endoscopic groups of $G$ imply part (2) for $G$. It was furthermore proved in \cite[Theorem 5.8]{KalGen} that for $p$-adic fields part (1') holds for quasi-split symplectic and orthogonal groups; the same argument also applies to quasi-split unitary groups.

\subsubsection{Automorphisms}
The next compatibility we discuss is with respect to automorphisms. Initially one may be interested in a particular connected reductive $F$-group $G'$ and wonder how an $F$\-automorphism $\theta'$ of $G'$ respects the (basic or refined) local Langlands correspondence. We will see however that this is a special case of the following more general consideration. 

Let as before $G$ be a quasi-split connected reductive $F$-group. Let $\theta \in \tx{Aut}_F(G)$. Then $\theta$ acts on $Z^1_\tx{bas}(\mc{E}^\tx{rig},G)$ via its action on $G$, and furthermore acts on the set of tuples $(G',\xi,z,\pi)$ by the rule $\theta(G',\xi,z,\pi) = (G',\xi\circ\theta^{-1},\theta(z),\pi)$. This induces an action of $\tx{Aut}_F(G)$ on the set of isomorphism classes of such tuples. The subgroup $\tx{Int}_F(G)$ of inner $F$-automorphisms of $G$ does \emph{not} act trivially, but it follows from Theorems \ref{thm:generic} and \ref{thm:central} and Remark \ref{rem:rx} that $\iota_{\varphi,\theta(\mf{w})}(\theta\dot\pi)=\iota_{\varphi,\mf{w}}(\dot\pi)$ for $\theta \in \tx{Int}_F(G)$ and $\dot\pi \in \Pi_\varphi$. Note this implies that $\theta$ preserves the compound $L$-packet $\Pi_\varphi$.

On the other hand the group $\tx{Aut}_F(\hat G)$ of $\Gamma$-equivariant automorphisms of $\hat G$ acts on the set of refined Langlands parameters by the rule $\theta(\varphi,\rho)=(\theta\circ\varphi,\rho\circ\theta^{-1})$, and the subgroup $\tx{Int}_F(\hat G)$ acts trivially on their $\hat G$-conjugacy classes by definition. Recall that the exact sequences of $\Gamma$-modules $1 \to \tx{Int}_{F^s}(G) \to \tx{Aut}_{F^s}(G) \to \tx{Out}_{F^s}(G) \to 1$ and $1 \to \tx{Int}(\hat G) \to \tx{Aut}(\hat G) \to \tx{Out}(\hat G) \to 1$ are split, and the choice of an $F$-pinning of $G$ resp. of a $\Gamma$-stable pinning of $\hat G$ determine $\Gamma$-equivariant splittings of these sequences. Recall finally that there is a natural identification $\tx{Out}(G) = \tx{Out}(\hat G)$, via which we obtain the identification $\tx{Aut}_F(G)/\tx{Int}_F(G)=\tx{Out}_F(G)=\tx{Out}_F(\hat G)=\tx{Aut}_F(\hat G)/\tx{Int}_F(\hat G)$.

To state the following conjecture, it is convenient to write $\iota(\mf{w},\dot\pi)=(\varphi,\rho)$, where $\varphi$ is the unique Langlands parameter with $\dot\pi \in \Pi_\varphi$, and $\rho = \iota_{\varphi,\mf{w}}(\dot\pi)$.

\begin{cnj} \label{cnj:auto}
Assume Conjecture \ref{cnj:refined}. Let $\varphi : L_F \to {^LG}$ be a Langlands parameter with compound $L$-packet $\Pi_\varphi$. For $\theta \in \tx{Out}_F(G) = \tx{Out}_F(\hat G)$,
\begin{enumerate}%
	\item $\Pi_{\theta\circ\varphi}=\{\theta\dot\pi\,|\, \dot\pi \in \Pi_\varphi\}$.
	\item $\iota(\theta(\mf{w}),\theta(\dot\pi)) = \theta(\iota(\mf{w},\dot\pi))$.
\end{enumerate}
\end{cnj}
Again part (1) can be formulated for $L$-packets on an individual group, and hence without assuming Conjecture \ref{cnj:refined}. In that form it appears as \cite[Lemma 6.18]{AV16} when $F=\R$.

Note the formal similarity between Conjectures \ref{cnj:contra} and \ref{cnj:auto}. If we let $\theta$ be the element of $\tx{Out}_F(G) = \tx{Out}_F(\hat G)$ that corresponds to the Chevalley involutions and set $\iota(\mf{w},\dot\pi)=(\varphi,\rho)$, then Conjecture \ref{cnj:contra}(2) states $\iota(\mf{w}^{-1},\dot\pi^\vee)=(\theta(\varphi),\theta(\rho)^\vee)$, while Conjecture \ref{cnj:auto}(2) states $\iota(\theta(\mf{w}),\theta(\dot\pi))=(\theta(\varphi),\theta(\rho))$. We are free to lift $\theta$ to an element of $\tx{Aut}_F(G)$ any way we like; if we take it to be the involution $\iota_{G,\mc{P}}$ defined by D. Prasad in \cite[Definition 1]{Dipendra19} with respect to a pinning related to the Whittaker datum $\mf{w}$ as in \cite[\S5.3]{KS99}, then $\theta(\mf{w})=\mf{w}^{-1}$. This shows that \cite[Conjecture 1]{Dipendra19} follows from Conjectures \ref{cnj:contra} and \ref{cnj:auto}.

Let us now discuss how Conjecture \ref{cnj:auto} gives information about the compatibility of the refined local Langlands correspondence with $F$-automorphisms of a fixed group $G'$. We realize $G'$ as a rigid inner twist $(\xi,z) : G \to G'$ of its quasi-split inner form $G$. Given $\theta' \in \tx{Aut}_F(G')$, the automorphism $\xi^{-1}\theta'\xi$ of $G$ need not respect the $F$-structure. However it does respect the $F^s$-structure and moreover the difference $(\xi^{-1}\theta'\xi)^{-1} \circ \sigma(\xi^{-1}\theta'\xi)$ is an inner automorphism of $G$ for each $\sigma \in \Gamma$. In other words, the image of $\xi^{-1}\theta'\xi$ in the group $\tx{Out}(G)=\tx{Aut}(G)/\tx{Int}(G)$ is an $F$-point. Let $\theta \in \tx{Aut}_F(G)$ be a lift of that $F$-point that preserves the Whittaker datum $\mf{w}$. Then $\xi^{-1}\theta'\xi = \tx{Ad}(g) \circ \theta$ for some $g \in G(F^s)$. Since the automorphism $\xi^{-1}\theta'\xi$ of $G$ commutes with the twisted action $\tx{Ad}(\bar z_\sigma) \sigma$ on $G(F^s)$ for all $\sigma \in \Gamma$, we see that the equality $\theta(\bar z_\sigma) = g^{-1}\bar z_\sigma \sigma(g)$ holds in $Z^1(F,G_\tx{ad})$. Then $y_e := (g^{-1} z_e \sigma_e(g)) \cdot \theta(z_e)^{-1} \in Z^1(\mc{E}^\tx{rig},Z(G))$ and we see that $(g,\theta')$ is an isomorphism $(G',\xi\circ\theta^{-1},y \cdot \theta(z),\pi) \to (G',\xi,z,\pi\circ\theta'^{-1})$. The class $[y]$ of $y$ is uniquely determined by $\theta'$ and $(\xi,z)$. From Theorem \ref{thm:central} we obtain

\begin{cor} \label{cor:auto}
Assume Conjectures \ref{cnj:refined} and \ref{cnj:auto}. Then
\[ \iota(G',\xi,z,\pi\circ\theta'^{-1}) = \theta(\iota(G',\xi,z,\pi)) \otimes (-d)^*([y]), \]
where the tensor product affects the second component $\rho$ of the refined parameter $(\varphi,\rho)$.
\end{cor}

The class $[y]$ is trivial if and only if $g$ can be chosen so that $(g,\theta')$ is an isomorphism $(G',\xi\circ\theta^{-1},\theta(z),\pi) \to (G',\xi,z,\pi\circ\theta'^{-1})$. This is not always possible: the simplest example is when $\theta'$ is an inner automorphism of $G'$ that comes from an element of $G'_\tx{ad}(F)$ which does not lift to $G'(F)$. In fact, this example is very useful and leads to the following application of Corollary \ref{cor:auto}.

\begin{cor} \label{cor:adjoint}
Assume Conjectures \ref{cnj:refined} and \ref{cnj:auto}. The action of $G'_\tx{ad}(F)$ on the set of representations of $G'(F)$ preserves each $L$-packet. More precisely, if $\theta'=\tx{Ad}(\bar g)$ for some $\bar g \in G'_\tx{ad}(F)$, then 
\[ \iota(G',\xi,z,\pi\circ\theta'^{-1})=\iota(G',\xi,z,\pi)\otimes (-d)^*([g^{-1}z_e\sigma_e(g)z_e^{-1}]), \]
for any lift $g \in G(F^s)$ of $\xi^{-1}(\bar g)$.
\end{cor}

\subsubsection{Homomorphism with abelian kernel and cokernel}

Let $f' : G_1' \to G_2'$ be a homomorphism of connected reductive $F$-groups with abelian kernel and cokernel, and $^Lf : {^LG_2'} \to {^LG_1'}$ be the corresponding $L$-homomorphism. For a Langlands parameter $\varphi_2 : L_F \to {^LG_2'}$ we can consider the composed parameter $\varphi_1 := {^Lf \circ \varphi_2}$ and the corresponding $L$-packets $\Pi_{\varphi_2}(G_2')$ and $\Pi_{\varphi_1}(G_1')$ provided by the basic local Langlands conjecture.

It is asserted in \cite[\S10]{BorCor} that for any $\pi_2 \in \Pi_{\varphi_2}(G_2')$ the representation $\pi_2 \circ f'$ of $G_1'(F)$ is a direct sum of finitely many members of $\Pi_{\varphi_1}(G_1')$. The refined local Langlands correspondence allows us to formulate a more precise expectation, namely about the multiplicity 
\[ m(\pi_1,\pi_2) = \tx{dim}\,\tx{Hom}(\pi_1,\pi_2\circ f') = \tx{dim}\,\tx{Hom}(\pi_2\circ f', \pi_1) \]
for each $\pi_1 \in \Pi_{\varphi_1}(G_1')$. To that end, note that $^Lf$ induces a map $\hat f : S_{\varphi_2}^+ \to S_{\varphi_1}^+$ via which we can define for any $\rho_1 \in \tx{Irr}(\pi_0(S_{\varphi_1}))$ and $\rho_2 \in \tx{Irr}(\pi_0(S_{\varphi_2}))$ the number
\[ m(\rho_1,\rho_2) = \tx{dim}\,\tx{Hom}(\rho_2,\rho_1\circ \hat f) = \tx{dim}\,\tx{Hom}(\rho_1\circ \hat f, \rho_2). \]
We can realize $G_i'$ as a rigid inner twist $(\xi_i,z_i) : G_i \to G_i'$ of its quasi-split inner form $G_i$ in such a way that $f = \xi_2^{-1} \circ f' \circ \xi_1$ is a homomorphism of $F$-groups $G_1 \to G_2$ and $z_2=f(z_1)$. 

\begin{cnj} \label{cnj:mult}
Let $\rho_i=\iota_{\varphi_i,\mf{w}}(G_i',\xi_i,z_i,\pi_i)$. Then
\[ m(\pi_1,\pi_2) = m(\rho_2,\rho_1). \]
\end{cnj}
This conjecture has been stated by Solleveld as \cite[Conjecture 2]{Soll20}, where it has been proven in some cases. A weaker form has been stated by Choiy in \cite{Choiy19}, where 
it has been proved under certain working hypotheses numbered 4.1, 4.3, 4.6, 4.11, and 4.15.  Among these, 4.1, 4.3, and 4.6 amount to Conjecture \ref{cnj:refined}, 4.11 is Corollary \ref{cor:adjoint}, and 4.15 is a certain dual version to 4.11 that also appears plausible. A direct verification of the basic version of this conjecture stated in \cite[\S10]{BorCor} in the setting of \cite{KalSLP} has been announced by Bourgeois--Mezo in \cite{BM21}, again assuming $F$ has characteristic zero.

\bibliographystyle{amsalpha}
\small
\bibliography{/Users/kaletha/Work/TexMain/bibliography.bib}

\providecommand{\bysame}{\leavevmode\hbox to3em{\hrulefill}\thinspace}
\providecommand{\MR}{\relax\ifhmode\unskip\space\fi MR }
\providecommand{\MRhref}[2]{%
  \href{http://www.ams.org/mathscinet-getitem?mr=#1}{#2}
}
\providecommand{\href}[2]{#2}
\begin{thebibliography}{CGH14}

\bibitem[ABV92]{ABV92}
Jeffrey Adams, Dan Barbasch, and David~A. Vogan, Jr., \emph{The {L}anglands
  classification and irreducible characters for real reductive groups},
  Progress in Mathematics, vol. 104, Birkh\"auser Boston, Inc., Boston, MA,
  1992. \MR{1162533 (93j:22001)}

\bibitem[Adl98]{Ad98}
Jeffrey~D. Adler, \emph{Refined anisotropic {$K$}-types and supercuspidal
  representations}, Pacific J. Math. \textbf{185} (1998), no.~1, 1--32.
  \MR{1653184 (2000f:22019)}

\bibitem[Art89a]{ArtIOR1}
James Arthur, \emph{Intertwining operators and residues. {I}. {W}eighted
  characters}, J. Funct. Anal. \textbf{84} (1989), no.~1, 19--84. \MR{999488
  (90j:22018)}

\bibitem[Art89b]{ArtUARC}
\bysame, \emph{Unipotent automorphic representations: conjectures},
  Ast\'erisque (1989), no.~171-172, 13--71, Orbites unipotentes et
  repr{\'e}sentations, II. \MR{1021499 (91f:22030)}

\bibitem[Art99]{Art99}
\bysame, \emph{On the transfer of distributions: weighted orbital integrals},
  Duke Math. J. \textbf{99} (1999), no.~2, 209--283. \MR{1708030 (2000i:22023)}

\bibitem[Art13]{Art13}
\bysame, \emph{The endoscopic classification of representations}, American
  Mathematical Society Colloquium Publications, vol.~61, American Mathematical
  Society, Providence, RI, 2013, Orthogonal and symplectic groups. \MR{3135650}

\bibitem[AS09]{AS09}
Jeffrey~D. Adler and Loren Spice, \emph{Supercuspidal characters of reductive
  {$p$}-adic groups}, Amer. J. Math. \textbf{131} (2009), no.~4, 1137--1210.
  \MR{2543925 (2011a:22018)}

\bibitem[AV92]{AV92}
Jeffrey Adams and David~A. Vogan, Jr., \emph{{$L$}-groups, projective
  representations, and the {L}anglands classification}, Amer. J. Math.
  \textbf{114} (1992), no.~1, 45--138. \MR{1147719 (93c:22021)}

\bibitem[AV16]{AV16}
\bysame, \emph{Contragredient representations and characterizing the local
  {L}anglands correspondence}, Amer. J. Math. \textbf{138} (2016), no.~3,
  657--682. \MR{3506381}

\bibitem[BDR17]{BDR17}
C\'edric Bonnaf\'e, Jean-Fran\c{c}ois Dat, and Rapha\"el Rouquier,
  \emph{Derived categories and {D}eligne-{L}usztig varieties {II}}, Ann. of
  Math. (2) \textbf{185} (2017), no.~2, 609--670. \MR{3612005}

\bibitem[BM21]{BM21}
Ad\`ele Bourgeois and Paul Mezo, \emph{Functoriality for supercuspidal
  {$L$}-packets}, 2021, preprint, arXiv:2109.09552.

\bibitem[BNS16]{BSN16}
A.~Bouthier, B.~C. Ng\^{o}, and Y.~Sakellaridis, \emph{On the formal arc space
  of a reductive monoid}, Amer. J. Math. \textbf{138} (2016), no.~1, 81--108.
  \MR{3462881}

\bibitem[Bor79]{BorCor}
A.~Borel, \emph{Automorphic {$L$}-functions}, Automorphic forms,
  representations and {$L$}-functions ({P}roc. {S}ympos. {P}ure {M}ath.,
  {O}regon {S}tate {U}niv., {C}orvallis, {O}re., 1977), {P}art 2, Proc. Sympos.
  Pure Math., XXXIII, Amer. Math. Soc., Providence, R.I., 1979, pp.~27--61.
  \MR{546608 (81m:10056)}

\bibitem[CGH14]{CGH14}
Raf Cluckers, Julia Gordon, and Immanuel Halupczok, \emph{Local integrability
  results in harmonic analysis on reductive groups in large positive
  characteristic}, Ann. Sci. \'{E}c. Norm. Sup\'{e}r. (4) \textbf{47} (2014),
  no.~6, 1163--1195. \MR{3297157}

\bibitem[Cho19]{Choiy19}
Kwangho Choiy, \emph{On multiplicity in restriction of tempered representations
  of {$p$}-adic groups}, Math. Z. \textbf{291} (2019), no.~1-2, 449--471.
  \MR{3936078}

\bibitem[CO21]{CO21}
Charlotte Chan and Masao Oi, \emph{Geometric {L}-packets of {H}owe-unramified
  toral supercuspidal representations}, preprint, arXiv:2105.06341, 2021.

\bibitem[Dil20]{Dillery20}
Peter Dillery, \emph{Rigid inner forms over local function fields}, 2020,
  preprint, arXiv:2008.04472.

\bibitem[Dil21]{Dillery21}
\bysame, \emph{Rigid inner forms over global function fields}, 2021, preprint,
  arXiv:2110.10820.

\bibitem[DL76]{DL76}
P.~Deligne and G.~Lusztig, \emph{Representations of reductive groups over
  finite fields}, Ann. of Math. (2) \textbf{103} (1976), no.~1, 103--161.
  \MR{0393266 (52 \#14076)}

\bibitem[DR09]{DR09}
Stephen DeBacker and Mark Reeder, \emph{Depth-zero supercuspidal {$L$}-packets
  and their stability}, Ann. of Math. (2) \textbf{169} (2009), no.~3, 795--901.
  \MR{2480618 (2010d:22023)}

\bibitem[DS18]{DS18}
Stephen DeBacker and Loren Spice, \emph{Stability of character sums for
  positive-depth, supercuspidal representations}, J. Reine Angew. Math.
  \textbf{742} (2018), 47--78. \MR{3849622}

\bibitem[Fin19]{Fin19}
Jessica Fintzen, \emph{On the construction of tame supercuspidal
  representations.}, preprint, arXiv:1908.09819, 2019.

\bibitem[Fin21]{Fin21}
\bysame, \emph{Types for tame {$p$}-adic groups}, Ann. of Math. (2)
  \textbf{193} (2021), no.~1, 303--346. \MR{4199732}

\bibitem[FKS21]{FKS}
Jessica Fintzen, Tasho Kaletha, and Loren Spice, \emph{A twisted {Y}u
  construction, {H}arish-{C}handra characters, and endoscopy}, arXiv:2106.09120
  (2021).

\bibitem[FOS20]{FOS20}
Yongqi Feng, Eric Opdam, and Maarten Solleveld, \emph{Supercuspidal unipotent
  representations: {L}-packets and formal degrees}, J. \'{E}c. polytech. Math.
  \textbf{7} (2020), 1133--1193. \MR{4167790}

\bibitem[FS21]{FS21}
Laurent Fargues and Peter Scholze, \emph{Geometrization of the local langlands
  correspondence}, preprint, arXiv:2102.13459, 2021.

\bibitem[G{\'e}r77]{gerardin:weil}
Paul G{\'e}rardin, \emph{Weil representations associated to finite fields}, J.
  Algebra \textbf{46} (1977), no.~1, 54--101. \MR{460477}

\bibitem[GHS]{GHS}
Wee~Teck Gan, Michael Harris, and Will Sawin, \emph{Local parameters of
  supercuspidal representations}, preprint, arXiv:2109.07737.

\bibitem[GL17]{GL17}
Alain Genestier and Vincent Lafforgue, \emph{Chtoucas restreints pour les
  groupes r\'eductifs et param\'etrisation de {L}anglands locale}, 2017,
  preprint, arXiv:1709.00978.

\bibitem[GL18]{GanLomeli18}
Wee~Teck Gan and Luis Lomel\'{\i}, \emph{Globalization of supercuspidal
  representations over function fields and applications}, J. Eur. Math. Soc.
  (JEMS) \textbf{20} (2018), no.~11, 2813--2858. \MR{3861809}

\bibitem[He08]{He08}
Xuhua He, \emph{On the affineness of {D}eligne-{L}usztig varieties}, J. Algebra
  \textbf{320} (2008), no.~3, 1207--1219. \MR{2427638 (2009c:20085)}

\bibitem[Hen93]{Hen93}
Guy Henniart, \emph{Caract\'{e}risation de la correspondance de {L}anglands
  locale par les facteurs {$\epsilon$} de paires}, Invent. Math. \textbf{113}
  (1993), no.~2, 339--350. \MR{1228128}

\bibitem[Hen00]{Hen00}
\bysame, \emph{Une preuve simple des conjectures de {L}anglands pour {$\mathrm{
  GL}(n)$} sur un corps {$p$}-adique}, Invent. Math. \textbf{139} (2000),
  no.~2, 439--455. \MR{1738446 (2001e:11052)}

\bibitem[Hen02]{Hen02}
\bysame, \emph{Une caract\'erisation de la correspondance de langlands locale
  pour $\mathrm{GL}(n)$}, Bulletin de la Soci\'et\'e Math\'ematique de France
  \textbf{130} (2002), no.~4, 587--602 (fr). \MR{1947454}

\bibitem[HII08]{HII08}
Kaoru Hiraga, Atsushi Ichino, and Tamotsu Ikeda, \emph{Correction to:
  ``{F}ormal degrees and adjoint {$\gamma$}-factors'' [{J}. {A}mer. {M}ath.
  {S}oc. \textbf{21} (2008), no. 1, 283--304; mr2350057]}, J. Amer. Math. Soc.
  \textbf{21} (2008), no.~4, 1211--1213. \MR{2425185}

\bibitem[HKT19]{HKT19}
Michael Harris, Chandrashekhar~B. Khare, and Jack~A. Thorne, \emph{A local
  langlands parameterization for generic supercuspidal representations of
  $p$-adic {$G_2$}}, 2019, preprint, arXiv:1909.05933.

\bibitem[HM08]{HM08}
Jeffrey Hakim and Fiona Murnaghan, \emph{Distinguished tame supercuspidal
  representations}, Int. Math. Res. Pap. IMRP (2008), no.~2, Art. ID rpn005,
  166. \MR{2431732 (2010a:22022)}

\bibitem[HT01]{HT01}
Michael Harris and Richard Taylor, \emph{The geometry and cohomology of some
  simple {S}himura varieties}, Annals of Mathematics Studies, vol. 151,
  Princeton University Press, Princeton, NJ, 2001, With an appendix by Vladimir
  G. Berkovich. \MR{1876802 (2002m:11050)}

\bibitem[Kal13]{KalGen}
Tasho Kaletha, \emph{Genericity and contragredience in the local {L}anglands
  correspondence}, Algebra Number Theory \textbf{7} (2013), no.~10, 2447--2474.
  \MR{3194648}

\bibitem[Kal15]{KalEpi}
\bysame, \emph{Epipelagic {$L$}-packets and rectifying characters}, Invent.
  Math. \textbf{202} (2015), no.~1, 1--89. \MR{3402796}

\bibitem[Kal16a]{KalSimons}
\bysame, \emph{The local {L}anglands conjectures for non-quasi-split groups},
  Families of automorphic forms and the trace formula, Simons Symp., Springer,
  2016, pp.~217--257. \MR{3675168}

\bibitem[Kal16b]{KalRI}
\bysame, \emph{Rigid inner forms of real and {$p$}-adic groups}, Ann. of Math.
  (2) \textbf{184} (2016), no.~2, 559--632. \MR{3548533}

\bibitem[Kal18a]{KalGRI}
\bysame, \emph{Global rigid inner forms and multiplicities of discrete
  automorphic representations}, Invent. Math. \textbf{213} (2018), no.~1,
  271--369. \MR{3815567}

\bibitem[Kal18b]{KalRIBG}
\bysame, \emph{Rigid inner forms vs isocrystals}, J. Eur. Math. Soc. (JEMS)
  \textbf{20} (2018), no.~1, 61--101. \MR{3743236}

\bibitem[Kal19a]{KalDC}
\bysame, \emph{On {$L$}-embeddings and double covers of tori over local
  fields}, arXiv:1907.05173 (2019).

\bibitem[Kal19b]{KalRSP}
\bysame, \emph{Regular supercuspidal representations}, J. Amer. Math. Soc.
  \textbf{32} (2019), no.~4, 1071--1170. \MR{4013740}

\bibitem[Kal19c]{KalSLP}
\bysame, \emph{Supercuspidal {$L$}-packets}, arXiv:1912.03274 (2019).

\bibitem[Kim07]{Kim07}
Ju-Lee Kim, \emph{Supercuspidal representations: an exhaustion theorem}, J.
  Amer. Math. Soc. \textbf{20} (2007), no.~2, 273--320. \MR{2276772
  (2008c:22014)}

\bibitem[Kot83]{Kot83}
Robert~E. Kottwitz, \emph{Sign changes in harmonic analysis on reductive
  groups}, Trans. Amer. Math. Soc. \textbf{278} (1983), no.~1, 289--297.
  \MR{697075 (84i:22012)}

\bibitem[Kot84]{Kot84}
\bysame, \emph{Stable trace formula: cuspidal tempered terms}, Duke Math. J.
  \textbf{51} (1984), no.~3, 611--650. \MR{757954 (85m:11080)}

\bibitem[Kot85]{Kot85}
\bysame, \emph{Isocrystals with additional structure}, Compositio Math.
  \textbf{56} (1985), no.~2, 201--220. \MR{809866 (87i:14040)}

\bibitem[Kot97]{Kot97}
\bysame, \emph{Isocrystals with additional structure. {II}}, Compositio Math.
  \textbf{109} (1997), no.~3, 255--339. \MR{1485921 (99e:20061)}

\bibitem[Kot14]{KotBG}
\bysame, \emph{{$B(G)$} for all local and global fields}, arXiv:1401.5728
  (2014).

\bibitem[Kre12]{Kret12}
Arno Kret, \emph{Existence of cuspidal representations of $p$-adic reductive
  groups}, 2012, preprint, arXiv:1205.2771.

\bibitem[KS99]{KS99}
Robert~E. Kottwitz and Diana Shelstad, \emph{Foundations of twisted endoscopy},
  Ast\'erisque (1999), no.~255, vi+190. \MR{1687096 (2000k:22024)}

\bibitem[KT18]{KT18}
Tasho Kaletha and Olivier Ta\"ibi, \emph{Global rigid inner forms vs
  isocrystals}, arXiv:1812.11373 (2018).

\bibitem[Lan89]{Lan89}
R.~P. Langlands, \emph{On the classification of irreducible representations of
  real algebraic groups}, Representation theory and harmonic analysis on
  semisimple {L}ie groups, Math. Surveys Monogr., vol.~31, Amer. Math. Soc.,
  Providence, RI, 1989, pp.~101--170. \MR{1011897 (91e:22017)}

\bibitem[Li19]{WenWeiLi19}
Wen-Wei Li, \emph{Contragredient representations over local fields of positive
  characteristic}, Algebra Number Theory \textbf{13} (2019), no.~5, 1197--1242.
  \MR{3981317}

\bibitem[LR87]{LR87}
R.~P. Langlands and M.~Rapoport, \emph{Shimuravariet\"aten und {G}erben}, J.
  Reine Angew. Math. \textbf{378} (1987), 113--220. \MR{895287 (88i:11036)}

\bibitem[LRS93]{LRS93}
G.~Laumon, M.~Rapoport, and U.~Stuhler, \emph{{$\mathscr{D}$}-elliptic sheaves
  and the {L}anglands correspondence}, Invent. Math. \textbf{113} (1993),
  no.~2, 217--338. \MR{1228127}

\bibitem[Lus95]{Lus95}
George Lusztig, \emph{Classification of unipotent representations of simple
  {$p$}-adic groups}, Internat. Math. Res. Notices (1995), no.~11, 517--589.
  \MR{1369407 (98b:22034)}

\bibitem[Lus02]{Lus02}
G.~Lusztig, \emph{Classification of unipotent representations of simple
  {$p$}-adic groups. {II}}, Represent. Theory \textbf{6} (2002), 243--289.
  \MR{1927955}

\bibitem[Mor88]{Mor88}
Lawrence Morris, \emph{{$P$}-cuspidal representations}, Proc. London Math. Soc.
  (3) \textbf{57} (1988), no.~2, 329--356. \MR{950594 (89j:22038)}

\bibitem[Mor89]{Mor89}
\bysame, \emph{{$P$}-cuspidal representations of level one}, Proc. London Math.
  Soc. (3) \textbf{58} (1989), no.~3, 550--558. \MR{988102 (90c:22056)}

\bibitem[MP94]{MP94}
Allen Moy and Gopal Prasad, \emph{Unrefined minimal {$K$}-types for {$p$}-adic
  groups}, Invent. Math. \textbf{116} (1994), no.~1-3, 393--408. \MR{1253198
  (95f:22023)}

\bibitem[MP96]{MP96}
\bysame, \emph{Jacquet functors and unrefined minimal {$K$}-types}, Comment.
  Math. Helv. \textbf{71} (1996), no.~1, 98--121. \MR{1371680 (97c:22021)}

\bibitem[MY20]{BMY20}
Alexander~Bertoloni Meli and Alex Youcis, \emph{An approach to the
  characterization of the local langlands correspondence}, 2020, preprint,
  arXiv:2003.11484.

\bibitem[Oha21]{Ohara21}
Kazuma Ohara, \emph{On the formal degree conjecture for non-singular
  supercuspidal representations}, 2021, preprint, arXiv:2106.00878.

\bibitem[Pra19]{Dipendra19}
Dipendra Prasad, \emph{Generalizing the {MVW} involution, and the
  contragredient}, Trans. Amer. Math. Soc. \textbf{372} (2019), no.~1,
  615--633. \MR{3968781}

\bibitem[Ree00]{Ree00}
Mark Reeder, \emph{Formal degrees and {$L$}-packets of unipotent discrete
  series representations of exceptional {$p$}-adic groups}, J. Reine Angew.
  Math. \textbf{520} (2000), 37--93, With an appendix by Frank L{\"u}beck.
  \MR{1748271 (2001k:22039)}

\bibitem[Ree08]{Ree08}
\bysame, \emph{Supercuspidal {$L$}-packets of positive depth and twisted
  {C}oxeter elements}, J. Reine Angew. Math. \textbf{620} (2008), 1--33.
  \MR{2427973 (2009e:22019)}

\bibitem[Ren10]{Ren10}
David Renard, \emph{Repr\'{e}sentations des groupes r\'{e}ductifs
  {$p$}-adiques}, Cours Sp\'{e}cialis\'{e}s, vol.~17, Soci\'{e}t\'{e}
  Math\'{e}matique de France, Paris, 2010. \MR{2567785}

\bibitem[Sch13]{Scholze13}
Peter Scholze, \emph{The local {L}anglands correspondence for {$\mathrm{GL}_n$}
  over {$p$}-adic fields}, Invent. Math. \textbf{192} (2013), no.~3, 663--715.
  \MR{3049932}

\bibitem[Sch21]{Schwein21}
David Schwein, \emph{Formal degree of regular supercuspidals}, 2021, preprint,
  arXiv:2101.00658.

\bibitem[Sha90]{Sha90}
Freydoon Shahidi, \emph{A proof of {L}anglands' conjecture on {P}lancherel
  measures; complementary series for {$p$}-adic groups}, Ann. of Math. (2)
  \textbf{132} (1990), no.~2, 273--330. \MR{1070599 (91m:11095)}

\bibitem[She82]{She82}
D.~Shelstad, \emph{{$L$}-indistinguishability for real groups}, Math. Ann.
  \textbf{259} (1982), no.~3, 385--430. \MR{661206 (84c:22017)}

\bibitem[She08a]{SheTE1}
\bysame, \emph{Tempered endoscopy for real groups. {I}. {G}eometric transfer
  with canonical factors}, Representation theory of real reductive {L}ie
  groups, Contemp. Math., vol. 472, Amer. Math. Soc., Providence, RI, 2008,
  pp.~215--246. \MR{2454336 (2011d:22013)}

\bibitem[She08b]{SheTE3}
\bysame, \emph{Tempered endoscopy for real groups. {III}. {I}nversion of
  transfer and {$L$}-packet structure}, Represent. Theory \textbf{12} (2008),
  369--402. \MR{2448289 (2010c:22016)}

\bibitem[She10]{SheTE2}
\bysame, \emph{Tempered endoscopy for real groups. {II}. {S}pectral transfer
  factors}, Automorphic forms and the {L}anglands program, Adv. Lect. Math.
  (ALM), vol.~9, Int. Press, Somerville, MA, 2010, pp.~236--276. \MR{2581952}

\bibitem[Sol19]{Soll19}
Maarten Solleveld, \emph{On unipotent representations of ramified $p$-adic
  groups}, 2019, preprint, arXiv:1912.08451.

\bibitem[Sol20]{Soll20}
Maarten Solleveld, \emph{Langlands parameters, functoriality and {H}ecke
  algebras}, Pacific J. Math. \textbf{304} (2020), no.~1, 209--302.
  \MR{4053201}

\bibitem[Spi18]{Spice18}
Loren Spice, \emph{Explicit asymptotic expansions for tame supercuspidal
  characters}, Compos. Math. \textbf{154} (2018), no.~11, 2305--2378.
  \MR{3867302}

\bibitem[Spi21]{Spice21}
\bysame, \emph{Explicit asymptotic expansions in p-adic harmonic analysis
  {II}}, preprint, arXiv:2108.12935, 2021.

\bibitem[SS13]{SS13}
Peter Scholze and Sug~Woo Shin, \emph{On the cohomology of compact unitary
  group {S}himura varieties at ramified split places}, J. Amer. Math. Soc.
  \textbf{26} (2013), no.~1, 261--294. \MR{2983012}

\bibitem[SZ18]{SZ18}
Allan Silberger and Ernst-Wilhelm Zink, \emph{{L}anglands classification for
  {L}-parameters}, Journal of Algebra \textbf{511} (2018), 299--357.

\bibitem[Ta{\"i}19]{TaibMult}
Olivier Ta{\"i}bi, \emph{Arthur's multiplicity formula for certain inner forms
  of special orthogonal and symplectic groups}, J. Eur. Math. Soc. (JEMS)
  \textbf{21} (2019), no.~3, 839--871. \MR{3908767}

\bibitem[Vog93]{Vog93}
David~A. Vogan, Jr., \emph{The local {L}anglands conjecture}, Representation
  theory of groups and algebras, Contemp. Math., vol. 145, Amer. Math. Soc.,
  Providence, RI, 1993, pp.~305--379. \MR{1216197 (94e:22031)}

\bibitem[Yu01]{Yu01}
Jiu-Kang Yu, \emph{Construction of tame supercuspidal representations}, J.
  Amer. Math. Soc. \textbf{14} (2001), no.~3, 579--622 (electronic).
  \MR{1824988 (2002f:22033)}

\end{thebibliography}

\end{document}